\documentclass[11pt,a4paper]{article}
\usepackage[utf8]{inputenc}
\usepackage{amsmath}
\usepackage{amsfonts}
\usepackage{amssymb}
\usepackage{mathtools}
\usepackage{amsthm}
\usepackage{enumitem}
\usepackage{xcolor}
\usepackage{hyperref}
\usepackage{verbatim}
\usepackage{mathrsfs}
\usepackage{titlefoot}

\newcommand{\set}[1]{\left\{#1\right\}}
\newcommand{\p}[1]{\left(#1\right)}
\newcommand{\va}[1]{\left|#1\right|}
\newcommand{\n}[1]{\left\|#1\right\|}

\newcommand{\C}{\mathbb{C}}
\newcommand{\N}{\mathbb{N}}
\newcommand{\R}{\mathbb{R}}

\DeclareMathOperator{\gl}{GL}
\DeclareMathOperator{\End}{End}

\newtheorem{thm}{Theorem}
\newtheorem{prop}{Proposition}[section]
\newtheorem{lemma}[prop]{Lemma}

\theoremstyle{definition}

\newtheorem{example}[prop]{Example}
\theoremstyle{remark}
\newtheorem{remark}[prop]{Remark}

\title{Rigorous justification for the space-split sensitivity algorithm to compute linear response in Anosov systems}
\date{}

\author{Nisha Chandramoorthy\footnote{Center for Computational Science and Engineering, Massachusetts Institute of Technology, Cambridge, MA 02139. Email: nishac@mit.edu}, Malo Jézéquel\footnote{Department of Mathematics, Massachusetts Institute of Technology, Cambridge, MA 02139. Email: mpjez@mit.edu}}

\begin{document}

\maketitle\unmarkedfntext{This is an author-created, un-copyedited version of an article accepted for publication in Nonlinearity.  The publisher is not responsible for any errors or omissions in this version of the manuscript or any version derived from it.  The Version of Record is available online at \url{https://doi.org/10.1088/1361-6544/ac7692}.}

\begin{abstract}
   Ruelle \cite{ruelleDifferentiationSRBStates1997,ruelleDifferentiationSRBStates2003,jiang_linear_response} gave a formula for linear response of transitive Anosov diffeomorphisms. Recently, practically computable realizations of Ruelle's formula have emerged that potentially enable sensitivity analysis of certain high-dimensional chaotic numerical simulations encountered in the applied sciences. In this paper, we provide full mathematical justification for the convergence of one such efficient computation, the space-split sensitivity, or S3, algorithm \cite{chandramoorthyEfficientComputationLinear2021}. In S3, Ruelle's formula is computed as a sum of two terms obtained by decomposing the perturbation vector field into a coboundary and a remainder that is parallel to the unstable direction. Such a decomposition results in a splitting of Ruelle's formula that is amenable to efficient computation. We prove the existence of the S3 decomposition and the convergence of the computations of both resulting components of Ruelle's formula. 
\end{abstract}

\section{Introduction}

\emph{Linear response} denotes the first order variation of statistical quantities associated to a dynamical system under a perturbation of this system. In the context of Anosov diffeomorphisms, this notion has been first studied by Ruelle, who gave a formula for linear response \cite{ruelleDifferentiationSRBStates1997, ruelleDifferentiationSRBStates2003,jiang_linear_response} (see \S \ref{sec:detailedIntro} for details).

The application of linear response theory has gained traction in the applied sciences, and especially in climate studies (see \cite{lucarini1,lucarini2,bodai,ghilLucarini} and references therein). Linear response is also interesting from the point of view of engineering, as it can be used to study design optimization, parameter selection and uncertainty quantification problems (see, e.g., \cite{huhn,fidkowski2011review,shimizu2018output,blonigan2014least,dow2011uncertainty} in the context of aerodynamic turbulent flows). But despite their importance, the lack of an efficient and rigorous computation of linear response has hindered downstream applications of linear response in high-dimensional chaotic models encountered in practice.

Recently, some new methods have been proposed to evaluate decomposed and regularized versions of Ruelle's formula \cite{chandramoorthyComputableRealizationRuelle2020,chandramoorthyEfficientComputationLinear2021, ni2021fast,sliwiak2021spacesplit}. The purpose of this paper is to give full mathematical justification of the validity of one of this method, the space-split sensitivity (S3) algorithm \cite{chandramoorthyEfficientComputationLinear2021,sliwiak2021spacesplit}. In doing so, we depart from the strategy of \cite{chandramoorthyEfficientComputationLinear2021}: the use of the spaces of anisotropic distributions from \cite{GLK} allows us to get a more direct understanding of the quantities that appear in the S3 algorithm. Such spaces have proved themselves extremely useful in hyperbolic dynamics since their introduction in \cite{BKL} (see \cite{quest} for a survey on this topic and \cite{Bal2} for a textbook presentation). Beyond the use of spaces of anisotropic distribution, our presentation differs from \cite{chandramoorthyEfficientComputationLinear2021} since we give a coordinate-free version of the ideas from \cite{chandramoorthyEfficientComputationLinear2021}.

Before recalling the steps of the S3 algorithm in \S \ref{sec:algorithm} and stating our main result, Theorem \ref{theorem:main}, let us recall a few facts about statistical properties of Anosov systems.

\subsection{Linear response of Anosov diffeomorphisms}
\label{sec:detailedIntro}

Let $\phi$ denote a $\mathcal{C}^r$ transitive Anosov diffeomorphism, with $r > 2$, on a compact $\mathcal{C}^\infty$ Riemannian manifold $M$. We recall that Anosov \cite[Definition 6.4.2]{katokhasselblatt} means that there is a splitting $TM = E_u \oplus E_s$ of the tangent bundle, invariant by the derivative of $\phi$, such that $E_u$ is uniformly expanded and $E_s$ is uniformly contracted by the derivative of $\phi$: there are constants $C > 0$ and $\lambda > 1$ such that for every $x \in M,n \in \mathbb{N}$ we have for all $v_u \in E_u(x)$ and $v_s \in E_s(x)$:
\begin{equation*}
\va{\mathrm{d}_x \phi^n \cdot v_s} \leq C \lambda^{-n} \va{v_s} \quad \textup{ and } \quad \va{\mathrm{d}_x \phi^{-n} \cdot v_u} \leq C \lambda^{-n} \va{v_u}.
\end{equation*}

The fundamental \emph{stable/unstable manifold theorem} \cite[\S 3.6]{yoccozIntroductionHyperbolicDynamics1995} asserts that the subbundle $E_u$ and $E_s$ are uniquely integrable: there are Hölder-continuous foliations $W_u$ and $W_s$ with uniformly $\mathcal{C}^r$ leaves (called respectively unstable and stable manifolds) that are tangent to $E_u$ and $E_s$ respectively.

Owing to their idealized setting, Anosov diffeomorphisms enjoy a neat statistical description. In particular, it is well-known  \cite{sinaiMarkovPartitionsDiffeomorphisms1968} that $\phi$ admits a unique SRB measure $\mu$, that is a unique invariant probability measure with smooth conditionals on unstable manifolds (see \cite{youngWhatAreSRB2002} for an introduction to and a survey of results on SRB measures). The SRB measure is {\em physical} in the sense that for Lebesgue almost every $x \in M$ and every continuous function $f$ on $M$,
\begin{equation}\label{eq:physical}
\frac{1}{n} \sum_{k = 0}^{n-1} f(\phi^k x)
\underset{n \to + \infty}{\to}
 \int_M f \mathrm{d}\mu.
\end{equation} 

Even though the measure $\mu$ is in general not absolutely continuous \cite[Corollary 1]{sinaiGIBBSMEASURESERGODIC}, its properties given above imply that $\mu$ governs the statistical properties of $\phi$ at the scale of Lebesgue measure, hence their importance in applied sciences. Many strong properties of the SRB measure $\mu$ are now well-known: central limit theorem \cite{denkerCentralLimitTheorem1989}, large deviation principle \cite{oreyDeviationsTrajectoryAverages2021}, full asymptotics for the correlations in term of resonances \cite{GLK}, etc. 

Of much importance for applications is the regularity of the dependence of $\mu$ on the dynamics $\phi$. To put it in more concrete terms, let $t \mapsto \phi_t$ be a smooth family of Anosov diffeomorphisms on $M$ with $\phi_0 = \phi$. For each $t$ near $0$, the Anosov diffeomorphism $\phi_t$ admits a unique SRB measure $\mu_t$ and it is then natural to wonder what is the regularity of the map $t \mapsto \mu_t$. When this map is differentiable at $0$, we say that \emph{linear response} holds. In that case, the value of the derivative at $0$ of the map $t \mapsto \mu_t$ is sometimes called the linear response, as it describes the first order response of statistical properties of a physical system under perturbations. 

In the context of transitive Anosov diffeomorphisms, linear response has been established by Ruelle \cite{ruelleDifferentiationSRBStates1997,ruelleDifferentiationSRBStates2003, jiang_linear_response}. In \cite{GLK}, Gou\"ezel and Liverani used the spectral properties of transfer operators to give a much more detailed description of the variation of the statistical properties of the SRB measure under perturbation. See also \cite[\S 5.3]{Bal2} for a textbook treatment of this problem. These references give a very detailed picture of the situation, let us just mention the following result that ensures that linear response holds.

\begin{thm}[{\cite[Theorem 2.8]{GLK}, \cite[Theorem 5.26]{Bal2}}]\label{theorem:linear_response}
Recall that $r > 2$. Let $t \mapsto \phi_t$ be a $\mathcal{C}^2$ function from a neighbourhood of $0$ in $\mathbb{R}$ to the space of $\mathcal{C}^r$ diffeomorhisms on the compact manifold $M$. Assume that $\phi = \phi_0$ is a transitive Anosov diffeomorphim. Then, for $t$ close enough to $0$ the map $\phi_t$ has a unique SRB measure $\mu_t$. 

Moreover, for every $\epsilon > 0$, there is $\eta > 0$ such that the map $t \mapsto \mu_t \in (\mathcal{C}^{2 + \epsilon})'$ is $\mathcal{C}^{1 + \eta}$ near $0$. In particular, if $f$ is a smooth function on $M$ then 
\begin{equation*}
t \mapsto \int_M f \mathrm{d}\mu_t
\end{equation*}
is $\mathcal{C}^{1+}$ near $0$.
\end{thm}

The above assumptions are rather restrictive, and cannot be verified with certainty in chaotic models encountered in practice. However, the application of linear response theory has been discussed for a variety of systems including climate models, turbulent flows and biological systems \cite{huhn, angxiuFluid, bodai, cessac}.  The justification for this robustness of linear response in noisy real-life systems has been studied from the statistical physics standpoint \cite{wormell1,wormell2}. Meanwhile, in the dynamical systems literature, rigorous proofs of existence of linear response have been found for a wider class of systems than Anosov systems, including certain partially hyperbolic systems \cite{dolgopyat}, some stochastic systems \cite{hairer,galatoloLinearResponseDynamical2017,bahsoun_ruziboev_saussol}, certain nonuniformly hyperbolic systems \cite{baladi-nonuniform} and even for some systems with intermittencies \cite{baladi-intermittent,bahsoun}. See the survey \cite{baladiLinearResponseElse2014} and references therein for further information on linear response.

The subject of our paper is the S3 algorithm \cite{chandramoorthyEfficientComputationLinear2021}, which computes the derivative at $0$ of the map $t \mapsto \mu_t$. To do so, we will rely on the following formula that has been proposed by Ruelle \cite{ruelleDifferentiationSRBStates1997,ruelleDifferentiationSRBStates2003, jiang_linear_response}: for a smooth observable $f$, we expect that
\begin{equation}\label{eq:Ruelle_formula}
\frac{\mathrm{d}}{\mathrm{d}t}\p{\int_M f \mathrm{d}\mu_t}_{|t= 0} = \sum_{n \geq 0} \mu\p{X\p{f \circ \phi^n}}.
\end{equation}
where the vector field $X$ describes the infinitesimal parameter perturbation:
\begin{equation*}
X \coloneqq \frac{\mathrm{d}}{\mathrm{d}t}\p{\phi_t}_{|t=0} \circ \phi^{-1},
\end{equation*}
so that the perturbed dynamics are given by $\phi_t = \phi + t X \circ \phi + {\cal O}(t^2)$. Notice that we used the common notation $\mu(f) = \int_M f \mathrm{d}\mu$. Under the assumption of Theorem \ref{theorem:linear_response}, Ruelle's formula \eqref{eq:Ruelle_formula} holds as soon as $f$ is $\mathcal{C}^{3 + \epsilon}$ (see for instance \cite[\S 5.3]{Bal2}).

At first glance, the series on the right-hand side of \eqref{eq:Ruelle_formula} does not seem to converge, due to the exponential growth of the derivatives of $\phi^n$ when $n$ tends to $+ \infty$. However, thanks to the integration against $\mu$, the series on the right-hand side is actually converging exponentially fast (see \cite{ruelleDifferentiationSRBStates1997,ruelleDifferentiationSRBStates2003,jiang_linear_response}; this can also be shown using the spectral investigation of a transfer operator associated with $\phi$). We will prove a slightly more general statement in Lemma \ref{lemma:convergence}. Nevertheless, it is inefficient to practically approximate the integrals against $\mu$ in \eqref{eq:Ruelle_formula} using a Monte Carlo method, i.e., approximating each integral as a finite sample average. The variance of the sample average also grows exponentially with $n$, making such a  direct evaluation of Ruelle's formula \eqref{eq:Ruelle_formula} computationally impractical for high-dimensional systems encountered in practice \cite{eyink, ensembleSensitivity}.
\subsection{The S3 algorithm}
\label{sec:algorithm}
The central object of this paper is the full mathematical justification for the use of the space-split sensitivity, or S3, algorithm proposed by Chandramoorthy and Wang in \cite{chandramoorthyComputableRealizationRuelle2020,chandramoorthyEfficientComputationLinear2021}. The S3 algorithm is one method to efficiently compute Ruelle's formula \eqref{eq:Ruelle_formula}. Notice that the S3 algorithm was first exposed in \cite{chandramoorthyEfficientComputationLinear2021} for system with one-dimensional unstable manifold and then extended to the general case by Sliwiak and Wang in \cite{sliwiak2021spacesplit}. We will also deal in this paper with the case of higher dimensional unstable direction. However, our approach to the issues specific to this more general case is not exactly the same as in \cite{sliwiak2021spacesplit}. Indeed, the methods to compute the unstable contribution (Step 3 of the algorithm below) in \cite{chandramoorthyEfficientComputationLinear2021,sliwiak2021spacesplit} are based on particular parametrization of the unstable manifolds, while we rather rely on the choice of a covariant derivative.

Ni's Linear Response Algorithm \cite{ni2021fast} is another alternative that converges to \eqref{eq:Ruelle_formula}, wherein the computations of the non-intrusive least squares shadowing algorithm \cite{angxiuNILSS} and a particular characterization of the unstable divergence are leveraged to provide a fast computation of \eqref{eq:Ruelle_formula}. Ni's Linear Response Algorithm proposes to split Ruelle's formula into shadowing and an unstable contribution while we propose a different decomposition in the S3 algorithm (see \cite[Appendix A]{chandramoorthyEfficientComputationLinear2021}, \cite{ni2021fast} for a more detailed comparison of these two algorithms). Abramov and Majda \cite{abramovMajda} developed a blended response technique for \eqref{eq:Ruelle_formula} in which they use a Gaussian approximation of the conditional densities of the SRB measures on the unstable manifolds to deal with terms corresponding to large values of $n$ in \eqref{eq:Ruelle_formula}.  

Other computational approaches to linear response are available, such as  \cite{bahsounRigorousComputationalApproach2017}, where, using suitable discretization schemes, the authors approximate the transfer operator associated to a one-dimensional dynamics and then use a functional analytic characterization of the linear response to approximate it. Our approach here is different, since, even if we use some knowledge on the spectrum of the transfer operator in \S \ref{section:convergence}, this is only needed for the justification of the convergence of the S3 algorithm and does not appear in the actual description of the algorithm below. We also introduce some transfer operators in \S \ref{section:spectral_radii}, but we only need to prove that their spectral radii are small to be able to solve certain equations (see Remark \ref{remark:solving_equations}), while the analysis from \cite{bahsounRigorousComputationalApproach2017} implies an understanding of the top eigenvalue of a transfer operator (including the existence of a spectral gap).

We recall that $r > 2$ and $\phi$ is a transitive $\mathcal{C}^r$ Anosov diffeomorphism on a compact $\mathcal{C}^\infty$ Riemannian manifold $M$ of dimension $d$. We shall denote by $d_u$ and $d_s$ respectively the unstable and stable dimension of $\phi$, that is the dimension of $E_u(x)$ and $E_s(x)$ respectively, for any $x \in M$. Considering the formula \eqref{eq:Ruelle_formula}, the S3 algorithm gives an approximation of the value of the sum
\begin{equation}\label{eq:result}
\Psi_{\phi}(X,f) \coloneqq \sum_{n \geq 0} \mu\p{X\p{f \circ \phi^n}}
\end{equation}
for any $\mathcal{C}^{1+}$ vector field $X$ and $\mathcal{C}^{1+}$ observable $f$. Here, $\mathcal{C}^{1+}$ just means $\mathcal{C}^{1+ \alpha}$ for some $\alpha > 0$.

The S3 algorithm is as follows. The input data are the map $\phi$, the vector field $X$ and the observable $f$. Further implementation details are discussed in Remark \ref{remark:implementation}.

\textbf{Step 1:} Find the unstable direction for $\phi$. To do so, start with any distribution of $d_u$-dimensional space $E_u^0$ and then let $E_u^{n+1}$ be the image of $E_u^0$ under the action of the derivative of $\phi^n$ for $n \geq 0$. Then, $E_u^n$ converges to $E_u$, as long as $E_u^0$ is chosen in an unstable cone for $\phi$. Notice that if we endow the Grassmanian of $d_u$-dimensional space over $M$ with a smooth Riemannian metric then the convergence is exponentially fast for the associated distance. For the practical computation, the $E_u^n$'s may be represented by one of its orthonormal basis, and then $E_u^{n+1}$ is deduced from $E_u^n$ by application of the derivative of $\phi$ followed by a Gram--Schmidt process. One must also compute ``the derivatives of $E_u$ in the unstable direction'' (see Remark \ref{remark:unstable_basis} for the meaning of this, and how to compute it).

\textbf{Step 2:} Compute a vector field $V$ so as to decompose $X$ as 
\begin{equation}\label{eq:decomposition_X}
X = \underbrace{Y}_{\textup{unstable term}} + \underbrace{V - \phi_*V}_{\textup{coboundary term}}, 
\end{equation}
where $Y$ is parallel to the unstable direction. Moreover, we require that $Y$ and $V$ are Hölder-continuous and that they have one Hölder-continuous derivative in the unstable direction (the precise meaning of this regularity assumption is detailed in \S \ref{section:spectral_radii}). The existence of the above decomposition is proved in Proposition \ref{prop:decomposition}, and a way to compute $Y$ and $V$ explicitly is explained in Remark \ref{remark:step_2}. The vector field $V$ is obtained by solving the linear equation $V - \mathcal{P} \phi_* V = X$ where $\mathcal{P}$ denotes the orthogonal projection on the orthogonal of the unstable direction for $\phi$. The lack of smoothness of $\mathcal{P}$ is the main reason we need to study objects that are smooth only in the unstable direction in \S \ref{section:spectral_radii}. 

We can then write
\begin{equation}\label{eq:decomposition_formula}
\Psi_{\phi}(X,f) = \Psi_{\phi}(Y,f) + \Psi_{\phi}(V - \phi_* V,f)
\end{equation}
and compute the contribution of $Y$ and $V - \phi_* V$ separately. Moreover, Lemma \ref{lemma:telescopic} asserts that $\Psi_{\phi}(V - \phi_* V,f) = \mu(Vf)$. The possibility to deal directly with $\Psi_{\phi}(V - \phi_* V,f)$ is one of the main novelty of our approach (the strategy from \cite{chandramoorthyEfficientComputationLinear2021} is a bit less direct). Notice that in \cite{chandramoorthyComputableRealizationRuelle2020,chandramoorthyEfficientComputationLinear2021,sliwiak2021spacesplit}, the coboundary contribution is called the ``stable contribution''. We do not use this convention here in order to highlight that there is no reason \emph{a priori} to have $V$ or $V - \phi_*V$ parallel to the stable direction.

\textbf{Step 3:} Find a Hölder-continuous function $\rho_Y$ such that for any $\mathcal{C}^{1 +}$ function $g$ we have $\mu(Yg) = \mu(g \rho_Y)$. The existence of $\rho_Y$ is obtained by integration by parts using the fact that $\mu$ has smooth conditionals on unstable manifolds (see Proposition \ref{prop:adjoint_Y}). In Lemma \ref{lemma:formula_rho}, we give a characterization of $\rho_Y$ that allows us to compute it explicitly (see also Remark \ref{remark:compute_rho}). Using $\rho_Y$, we can obtain the unstable contribution as
\begin{equation}\label{eq:unstable_contribution}
\Psi_{\phi}(Y,f) = \sum_{n \geq 0} \mu\p{f \circ \phi^n \rho_Y},
\end{equation}
where the sum converges exponentially fast due to the exponential decay of correlations (see for instance \cite[Theorem 7.11]{Bal2}) for $\phi$ (notice that $\rho_Y$ is Hölder and that $\mu(\rho_Y) = \mu(Y1) = 0$). The advantage of \eqref{eq:unstable_contribution} over the initial definition of $\Psi_{\phi}(Y,f)$ is that we do not need to compute exponentially growing derivatives anymore. Due to the exponential convergence of the series in \eqref{eq:unstable_contribution}, we get a good approximation of $\Psi_{\phi}(Y,f)$ by truncating the series at some $L \gg 1$:
\begin{equation}\label{eq:approximation_unstable}
\Psi_{\phi}(Y,f) \simeq \sum_{n = 0}^{L - 1} \mu\p{f \circ \phi^n \rho_Y}.
\end{equation}

\textbf{Step 4:} Compute $\mu(f \rho_Y),\dots, \mu( f \circ \phi^{L-1} \rho_Y)$ and $\mu(Vf)$ using either of the two methods (probabilistic or deterministic) proposed in \S \ref{section:monte_carlo}, and take as an approximate value for $\Psi_{\phi}(X,f)$:
\begin{equation}\label{eq:approximation_S3}
\Psi_{\phi}(X,f) \simeq \mu(Vf) + \sum_{n =0}^{L-1} \mu\p{ f \circ \phi^n \rho_Y}.
\end{equation}

Our main result is then the following.

\begin{thm}\label{theorem:main}
The S3 algorithm converges: with the notation above, we have
\begin{equation}\label{eq:S3_formula}
\Psi_{\phi}(X,f) = \mu(Vf) + \sum_{n \geq 0} \mu(f \circ \phi^n \rho_Y),
\end{equation}
where we recall that $X$ and $f$ are $C^{1+}$, while $\phi$ is $\mathcal{C}^{2+}$.
\end{thm}

The identity \eqref{eq:S3_formula} is an immediate consequence of Propositions \ref{prop:decomposition} and \ref{prop:adjoint_Y} and Lemma \ref{lemma:convergence}. The speed of convergence of the S3 algorithm is discussed in \S \ref{section:monte_carlo}. 
\begin{remark}
The S3 algorithm described above is the same as in \cite{chandramoorthyEfficientComputationLinear2021}. However, the decomposition \eqref{eq:decomposition_X} is not highlighted in the latter (while it is already present, see equation (8.23) there), because the authors could not prove directly the convergence of the term $\Psi_\phi(V - \phi_* V, f)$ in the right hand side of \eqref{eq:decomposition_formula}. Hence, \cite{chandramoorthyEfficientComputationLinear2021} expresses \eqref{eq:decomposition_formula} (see (4.2) in \cite{chandramoorthyEfficientComputationLinear2021}) as the asymptotic limit (in the length of the trajectory) of a decomposition that emerges from the same algorithm. However, in effect, the linear response $\Psi_\phi(X,f)$ is already decomposed in \cite{chandramoorthyEfficientComputationLinear2021} as the sum of an unstable contribution, denoted by ``$\langle J,\partial_s \mu_s \rangle^u$'', and a stable contribution (which we call coboundary contribution), denoted by ``$\langle J,\partial_s \mu_s \rangle^s$'' (see equation (4.2) there). The stable contribution is shown to coincide with the term $\mu(Vf)$ in \eqref{eq:S3_formula} (see equation (5.4) in \cite{chandramoorthyEfficientComputationLinear2021}), and it is stated that the unstable contribution coincides with the sum in \eqref{eq:S3_formula} (see equation (6.18) in \cite{chandramoorthyEfficientComputationLinear2021}). Notice however that there is a gap (in \S 6.4 of \cite{chandramoorthyEfficientComputationLinear2021}) in the proof of the latter statement in \cite{chandramoorthyEfficientComputationLinear2021}. The estimates that we establish in \S \ref{section:spectral_radii} could be used to fill that gap.

The main novelty of our work is the use of the spaces of anisotropic distributions from \cite{GLK} that allows us to prove the convergence of $\Psi_\phi(V-\phi_* V,f)$ (see Lemma \ref{lemma:telescopic}) and thus to work directly with the decomposition \eqref{eq:decomposition_X}. We provide consequently a more direct analysis of each term in \eqref{eq:S3_formula}. 

Moreover, we give in \S \ref{section:spectral_radii} a detailed discussion of classes of regularity tailored for our problem: we study objects that have ``Hölder derivatives in the unstable direction'' without \emph{a priori} being differentiable. Some functions that are differentiable only in the unstable direction already appeared in \cite{chandramoorthyEfficientComputationLinear2021}, but we give here a more systematic and careful analysis of these issues. This allows us to isolate a general mechanism which is behind the definition and computation of most objects that appear in the S3 algorithm (see Proposition \ref{prop:bound_spectral_radius} and Remark \ref{remark:solving_equations}).

Our hope is that our approach clarifies the regularity issues involved in the S3 algorithm. Notice also that our analysis is coordinate-free, and tries to unveil the intrinsic notion that backs the validity of the S3 algorithm.
\end{remark}

\begin{remark}\label{remark:implementation}
Let us describe briefly how the S3 algorithm can be actually implemented (see \cite{chandramoorthyEfficientComputationLinear2021,sliwiak2021spacesplit} for details). In order to approximate the integrals $\mu(f \rho),\dots, \mu( f \circ \phi^{L-1} \rho)$ and $\mu(Vf)$ in Step 4, we only need to know $\rho$ and $V$ at a finite number of points. To compute $\rho$ and $V$ at a point $x$, we can proceed as follows.

We start by computing a large segment $\phi^{-n}x,\dots,\phi^{-1}x,x$ of the backward orbit of $x$. Then, we cut this orbit in $3$ parts of comparable length $O_1 = \set{\phi^{-n}x,\dots, \phi^{- m-1}x}$, $O_2 = \set{\phi^{-m}x,\dots,\phi^{- \ell - 1} x}$ and $O_3 = \set{\phi^{-\ell} x,\dots,x}$. We use the first orbit segment $O_1$ to perform Step 1: we choose any $d_u$-dimensional plane $E_u^0(\phi^{-n}x)$ over $\phi^{-n}x$ (if possible near the unstable direction), we can then recursively compute $E_u^1(\phi^{-n+1}x),\dots, E_u^{k}(\phi^{-n+ k}x)$. If $O_1$ is large enough, this process gives us a good approximation of $E_u$ at the points of $O_2$ and $O_3$ (the error is exponentially small in the size of $O_1$, and the convergence holds provided the initial direction does not intersect the stable direction).

In Step 1, one also needs to compute ``the derivatives of $E_u$ in the unstable direction'' (see Remark \ref{remark:unstable_basis} for the meaning of this and how to compute it). It is necessary to know $E_u$ to perform this computation, but it is not needed until Step 3, so that one can use the orbit segment $O_2$ to compute it at the point of $O_3$. Since this is a computation based on the method from Remark \ref{remark:solving_equations}, the error made is exponential in the length of the orbit segment used.

We use the orbit segment $O_2$ to perform Step 2. The method given in Remark \ref{remark:step_2} (see also Remark \ref{remark:solving_equations}) to compute $Y$ and $V$ at a point $y \in M$ only uses values of $X,E_u$ and of the derivative of $\phi$ at points in the backward orbit of $y$. This procedure converges to $V$ and $Y$ exponentially fast in the length of the orbit used. Hence, using the points of $O_2$ (with the value of $E_u$ computed at Step 1), we get an approximation of $Y$ and $V$ at the points of $O_3$.

Finally, we proceed at Step 3 using the orbit segment $O_3$: using a similar method as in Step 2, we can compute $\rho_Y$ at $x$ using the values of $Y$ and $E_u$ at the points of $O_3$ (see Remark \ref{remark:compute_rho}), and the error made is exponentially small in the length of $O_3$.

Now that we know how to compute the values of $V$ and $\rho_Y$ at a given point, we can apply the methods described in \S \ref{section:monte_carlo} to approximate the integrals $\mu(f \rho_Y),\dots, \mu( f \circ \phi^{L-1} \rho_Y)$ and $\mu(Vf)$. In practice, the probabilistic method in \S \ref{subsection:Montecarlo} do not need the values of $V$ and $\rho_Y$ at any specific $x$, but rather along any $\mu$-typical orbit of $\phi$. The deterministic method in \S \ref{subsection:deter} only requires to compute these functions at $\phi^n x$ for some $x$'s and $n$ large. As a result, one need not invert the map $\phi$, but only compute forward orbits (and use the forward orbit of $x$ as a backward orbit for $\phi^n x$). This is essential when working with a hyperbolic attractor (see \S \ref{section:hyperbolic_attractors}).
\end{remark}

\begin{remark}
Instead of \eqref{eq:decomposition_X}, one could naturally imagine decomposing $X$ as $X = V_u + V_s$ where $V_u$ and $V_s$ are respectively in the unstable and stable direction. Then, the terms in the sum \eqref{eq:result} defining $\Psi_\phi(V_s,f)$ would be exponentially decaying. Further, one could use integration by parts on the unstable manifold to express $\Psi_\phi(V_u,f)$ as an exponentially decaying sum of correlations involving an unstable divergence denoted by ${\rm div}^u V_u$  in \cite{ruelleDifferentiationSRBStates1997}. This divergence with respect to the volume on the unstable manifold defined by the SRB measure has been shown to be H\"older continuous (see \cite{jiang_linear_response}), even though the derivatives of $V_u$ in the unstable direction do not exist. Our approach  in the S3 algorithm rests on a recursive, convergent procedure to differentiate functions in the unstable direction, but the application of this procedure to compute the unstable divergence is precluded by the latter fact. Hence, the algorithm cannot be used straightforwardly with the decomposition into $V_u + V_s.$ It is worth noting that Ni's fast linear response algorithm \cite{ni2021fast} ends up computing the unstable divergence term, but approximating instead an alternative expression for the divergence that does not explicitly involve differentiating $V_u$ in the unstable direction. Ni's \cite{ni2021fast} decomposition is into a ``shadowing'' direction and an unstable vector field. This decomposition leads, through integration by parts on the unstable manifold, to the problematic unstable divergence term, which we bypass here thanks to a decomposition \eqref{eq:decomposition_formula} that is differentiable in the unstable directions. 
\end{remark}

The remainder of this article proves the main result (Theorem \ref{theorem:main}). The organization is as follows. In \S \ref{section:spectral_radii}, we give a precise meaning to the expression ``having a Hölder-continuous derivative in the unstable direction''. We also give a bound on the spectral radii of certain transfer operators, Proposition \ref{prop:bound_spectral_radius}, that will be the main tool in order to construct $Y,V$ and $\rho$ in Steps 2 and 3 of the S3 algorithm. This bound is also what ensures that the $Y,V$ and $\rho$ may be computed with exponential precision. In \S \ref{section:Y_and_Z}, we explain how the vector fields $Y$ and $V$ from Step 2 are constructed (Proposition \ref{prop:decomposition}). The section \S \ref{section:unstable} is devoted to the existence, Proposition \ref{prop:adjoint_Y}, and computation, Lemma \ref{lemma:formula_rho}, of the function $\rho$ from Step 3. In \S \ref{section:convergence}, we prove the formula $\Psi_{\phi}(V - T_* V,f) = \mu(Vf)$ (see Lemma \ref{lemma:telescopic}). This is one of the main novelty from our work, which makes possible a more direct approach to the S3 algorithm decomposition. The proof of this crucial result is based on the spectral decomposition of a Koopman operator acting on the spaces of anisotropic distributions from \cite{GLK}. In \S \ref{section:monte_carlo}, we explain how the integrals in Step 4 of the algorithm may be approximated, and we estimate the precision of the S3 algorithm (see Remarks \ref{remark:precision_monte_carlo} and \ref{remark:precision_Riemann}). Finally, in \S \ref{section:hyperbolic_attractors}, we sketch how to generalize our proof of the validity of the S3 algorithm to the case of hyperbolic attractors.

\subsection*{Acknowledgements} Both authors thank Semyon Dyatlov for introducing them to each other and for useful discussions about this work. N.C. is supported by the Air Force Office of Scientific Research Grant No. FA8650-19-C-2207 and Department of Energy Research Grant No. DE-FOA-0002068-0018. N.C. gratefully acknowledges Qiqi Wang for the S3 algorithm in \cite{chandramoorthyEfficientComputationLinear2021} and numerous valuable inputs. During the beginning of the elaboration of this work, M.J. was supported by the European Research Council (ERC) under the European Union's Horizon 2020 research and innovation programme (grant agreement No 787304). Part of this work was done while MJ was working at LPSM\footnote{Laboratoire de Probabilit\'es, Statistique et Mod\'elisation (LPSM), CNRS, Sorbonne Universit\'e, Universit\'e de Paris, 4, Place Jussieu, 75005 Paris, France}.

\section{Spectral radii of transfer operators on vector bundles}\label{section:spectral_radii}

In this section, we define precisely what we mean by having a Hölder-continuous derivative in the unstable direction (as in Step 2 of the S3 algorithm) and give the basic properties of this regularity condition (see Lemma \ref{lemma:regularity}). The main point of this section is Proposition \ref{prop:bound_spectral_radius} that gives a bound on the spectral radii of certain transfer operators associated with $\phi$. This result extracts and generalizes the main mechanism underlying in the proofs of Lemmas 8.6,8.7,8.8 and 8.11 and Proposition 8.9 in \cite{chandramoorthyEfficientComputationLinear2021}. Since we will need to consider transfer operators acting both on functions and on vector fields, it is convenient to work here in the context of a general vector bundle. Considering the action of transfer operators on vector bundle is standard, see for instance \cite[\S 6.4]{Bal2}.

Let $p : E \mapsto M$ be a $\mathcal{C}^r$ complex vector bundle over $M$. We denote by $E_x \coloneqq p^{-1}(\set{x})$ the fiber of $E$ over $x$. Let $\Phi$ be a lift of $\phi$, that is a map from $E$ to itself that satisfies $p \circ \Phi = \phi \circ p$, and whose restriction to each fiber is linear. To $\Phi$, we associate a transfer operator acting on sections of $E$, defined by
\begin{equation}\label{eq:general_transfer_operator}
\mathcal{L}_{\Phi} v (x) = \Phi\p{v (\phi^{-1}x)},
\end{equation}
for $x \in M$ and $v$ a section of $E$. Before studying the properties of operators of the form $\mathcal{L}_{\Phi}$, let us give the main examples of such operators that appear in the S3 algorithm.

\begin{example}\label{example:scalar}
If $E = M \times \C$ is the trivial line bundle over $M$ and $h : M \to \C$ is a scalar function, then we can define a lift $\Phi$ for $\phi$ by
\begin{equation*}
\Phi(x,z) = \p{\phi x , h(\phi x) z} \quad \textup{for} \quad (x,z) \in M \times \C
\end{equation*}
In this case, the sections of $E$ identify with the functions from $M$ to $\mathbb{C}$, and the operator $\mathcal{L}_{\Phi}$ is given by the formula
\begin{equation*}
\mathcal{L}_\Phi u(x) = h(x) u(\phi^{-1}x),
\end{equation*}
for $x \in M$ and $u$ a function on $M$.
\end{example}

\begin{example}\label{example:vector_fields}
Another basic example is $E = T M$, the tangent vector bundle of $M$, we can then define the lift $\mathrm{d}\phi$ of $\phi$ by 
\begin{equation*}
\Phi(x,v) = \mathrm{d}_x \phi \cdot v \quad \textup{ for } \quad v \in T_x M.
\end{equation*}
This example will be crucial in the decomposition of $X$. Indeed, in this case we retrieve $\mathcal{L}_{\Phi} = \phi_*$, defined by,
\begin{equation*}
\mathcal{L}_\Phi V(x) = \phi_* V(x) = \mathrm{d}_{\phi^{-1}x} \phi \cdot
V(\phi^{-1} x),
\end{equation*}
for $x \in M$ and $V$ a vector field on $M$.
\end{example}

\begin{example}
From the lift $\mathrm{d}\phi$ of $\phi$ in Example \ref{example:vector_fields}, one can construct the lift $\Phi = \mathcal{P}\mathrm{d}\phi$ where, for each $x \in M$, the operator $\mathcal{P}$ acts on $T_x M$ like the orthogonal projection on the orthogonal complement of the unstable direction $E_u(x)$. For $V$ a vector field on $M$, we see that $\mathcal{L}_\Phi V$ is obtained by projecting $\phi_* V$ orthogonally  on the orthogonal complement of $E_u$ in each fiber. This lift is central in the S3 algorithm, see the proof of Proposition \ref{prop:decomposition} below. The map $\mathcal{P}$ is a basic example of object that is not smooth, but has derivatives in the unstable direction (see Lemma \ref{lemma:regularity}).
\end{example}

\begin{remark}
Notice that if $\Phi : E \to E$ is a lift of $\phi$, then the composition $\Phi^m : E \to E$ is a lift of $\phi^m$, and we have $\mathcal{L}_{\Phi^m} = (\mathcal{L}_{\Phi})^m$.
\end{remark}

\begin{remark}\label{remark:operation_bundle}
If $E$ is a vector bundle over $M$, we will denote by $E^*$ the dual vector bundle, i.e. the natural vector bundle such that, for every $x \in M$, the fiber of $E^*$ over $x$ is the dual of the fiber of $E$ over $x$. If $F$ is another vector bundle, we recall that the fiber over $x$ of the vector bundles $E \oplus F$ and $E \otimes F$ are respectively $E_x \oplus F_x$ and $E_x \otimes F_x$.

Recall that if $E_1,\dots,E_N,F$ are vector bundle over $M$, and if $v_1,\dots,v_N$ are section of $E_1,\dots,E_N$ respectively and $\gamma$ is a section of the product $E_1^* \otimes \dots \otimes E_N^* \otimes F$, then we can define a section $\gamma(v_1,\dots,v_N)$ of $F$, since for every $x \in M$, the vector $\gamma(x)$ identifies with a $N$-linear form from $E_{1,x} \times \dots \times E_{N,x}$ to $F_x$.

If $E$ is a vector bundle over $M$, then we will denote by $\phi^* E$ the pullback of $E$ by $\phi$, this is a vector bundle such that, for every $x \in M$, the fiber of $\phi^* E$ over $x$ is $E_{\phi x}$. Notice that, if $v$ is a section of $E$, then $v \circ \phi$ is a section of $\phi^* E$. We define similarly $(\phi^{-1})^*E$.

For a more detailed discussion of these standard constructions, the reader may refer for instance to \cite[\S 1.F]{ghl_geometry}, \cite[Chapter 3]{husemoller}.
\end{remark}

Let us fix a $\mathcal{C}^r$ covariant derivative $\nabla$ on $E$. By this, we just mean that if $X$ is a $\mathcal{C}^r$ vector field on $M$ and $v$ a $\mathcal{C}^r$ section of $E$ then $\nabla_{X} v$ is a $\mathcal{C}^{r-1}$ section of $E$, and we require in addition that if $f$ is a $\mathcal{C}^r$ function on $M$ then $\nabla_{fX}v = f \nabla_X v$ and $\nabla_X (fv) = Xf.v + f \nabla_Xv$ (where $Xf$ denotes the Lie derivative of $X$ applied to $f$, that is $Xf(x) = \mathrm{d}_x f \cdot X(x)$ for $x \in 
M$), and that $\nabla$ is $\mathbb{R}$-bilinear. Notice that if we work in a trivialization of $E$ that identifies locally the fibers of $E$ with $\R^m$, with canonical basis $e_1,\dots,e_m$, and the tangent space of $M$ with $\R^n$, with canonical basis $\partial_{x_1},\dots,\partial_{x_n}$, then, writing $X = \sum_{j = 1}^n X^j \partial_{x_j}$ and $v = \sum_{k = 1}^m v^k e_k$, we have
\begin{equation}\label{eq:expression_coordinate}
\nabla_{X} v = \mathrm{d}v \cdot X + \sum_{\substack{i,k = 1,\dots,m \\ j = 1,\dots,n}} \Gamma_{j,k}^i X^j v^k e_i, 
\end{equation}
where the Christoffel symbols $\Gamma_{j,k}^i$ are $\mathcal{C}^{r-1}$ functions defined by the relation
\begin{equation*}
\nabla_{\partial_{x_j}}(e_k) = \sum_{i=1}^m \Gamma_{j,k}^i e_i.
\end{equation*}
Here the differential $\mathrm{d}v$ makes sense since the trivialization of $E$ identifies $v$ with a function from an open subset of $\mathbb{R}^n$ to $\mathbb{R}^m$.  If $E$ is as in Example \ref{example:scalar}, we can just take $\nabla$ to be the trivial connection (given by $\nabla_X f = X f$), and in Example \ref{example:vector_fields} it is natural to choose the Levi--Civita connection associated to the Riemannian metric on $M$ (see for instance \cite[Theorem 2.51 and Definition 2.53]{ghl_geometry}).

\begin{remark}\label{remark:different_connection}
Notice that if $\nabla$ is a $\mathcal{C}^r$ covariant derivative on $E$, then there is a (unique) $\mathcal{C}^r$ covariant derivative $\nabla$ on $E^*$, the bundle whose fibers are the duals of the fibers of $E$, such that if $l$ and $v$ are sections respectively of $E^*$ and $E$ then we have
\begin{equation*}
X(l(v)) = (\nabla_X l)(v) + l(\nabla_X v)
\end{equation*}
for any $\mathcal{C}^r$ vector field $v$. The uniqueness may be deduced from the expression $(\nabla_X l)(v) = X(l(v)) - l(\nabla_X v)$ and the existence only needs to be checked locally, which can be done using the expression \eqref{eq:expression_coordinate}. Using the dual basis of $e_1,\dots,e_m$ to trivialize $E^*$, the Christoffel symbols $\widetilde{\Gamma}_{j,k}^i$ for $E^*$ are deduced from those of $E$ by the formula $\widetilde{\Gamma}_{j,k}^i = - \Gamma_{j,i}^k$.

Notice that if $E = TM$ is the tangent vector bundle of $M$, endowed with the Levi--Civita covariant derivative \cite[Theorem 2.51 and Definition 2.53]{ghl_geometry}, and $E^* = T^*M$ is the cotangent vector bundle, then the connection on $T^* M$ obtained by the procedure above is the same as the Levi--Civita connection under the identification $TM \simeq T^*M$ given by the Riemannian metric on $M$.

Similarly, if $F$ is another $\mathcal{C}^r$ vector bundle, associated with a $\mathcal{C}^r$ covariant derivative that we also denote by $\nabla$, then there is a unique $\mathcal{C}^r$ covariant derivative on the vector bundle $E \otimes F$ such that
\begin{equation*}
\nabla_X (u \otimes v) = (\nabla_X u) \otimes v + u \otimes (\nabla_X v),
\end{equation*}
for any $\mathcal{C}^r$ sections $u$ and $v$ respectively of $E$ and $F$, and any $\mathcal{C}^r$ vector field $X$. We can also define a covariant derivative on $E \oplus F$ such that $\nabla_X(u \oplus v) = \nabla_X u \oplus \nabla_X v$.

Identifying the bundle $\End(E)$, whose fibers are made of the endomorphisms of the fibers of $E$, with $E \otimes E^*$, we get a $\mathcal{C}^r$ covariant derivative on $\End(E)$ that satisfies the useful relation
\begin{equation}\label{eq:Leibniz_endomorphism}
\nabla_X(AB) = (\nabla_X A)B + A (\nabla_X B),
\end{equation}
for any $\mathcal{C}^r$ sections $A$ and $B$ of $\End(E)$ and any $\mathcal{C}^r$ vector field $X$ on $M$.

More generally, if $\gamma$ and $v_1,\dots,v_N$ are as in Remark \ref{remark:operation_bundle} and $X$ is a vector bundle, then we have the Leibniz relation:
\begin{equation}\label{eq:Leibniz}
\begin{split}
\nabla_{X}\p{\gamma(v_1,\dots,v_N)} & = \nabla_X \gamma(v_1,\dots,v_N) \\ & \qquad \qquad + \sum_{j = 1}^N \gamma(v_1,\dots,v_{j-1}, \nabla_X v_j, v_{j+1},\dots,v_N).
\end{split}
\end{equation}

The reader may refer to \cite[\S 2.B.3]{ghl_geometry} for a discussion of covariant derivatives on tensors. Finally, notice that one may define a connection on $(\phi^{-1})^* E$ by setting $\nabla_X v = \p{\nabla_{\phi^{-1}_* X} (v \circ \phi)} \circ \phi^{-1}$ for a section $v$ of $(\phi^{-1})^* E$ and a vector field $X$. 
\end{remark}

We will study the operator $\mathcal{L}_{\Phi}$ acting on spaces that we introduce now. The first space is the space $\mathcal{C}^0\p{M;E}$ of continuous sections of $E$. In order to define the topology on $\mathcal{C}^0\p{M;E}$, we choose a smooth Riemannian metric on $E$ and, letting $\va{\cdot}$ denote the associated norm on the fibers of $E$, we define the norm on $\mathcal{C}^0\p{M;E}$ by
\begin{equation*}
\va{v}_{\infty} = \sup_{x \in M} \va{v(x)} \quad \textup{for} \quad v \in \mathcal{C}^0\p{M;E}.
\end{equation*}
One can then check that the norm $\va{\cdot}_\infty$ induces a structure of Banach space on $\mathcal{C}^0\p{M;E}$. Moreover, another choice of Riemannian metric on $E$ leads to an equivalent norm.

We will also need for $0 < \alpha < 1$ the space $\mathcal{C}^\alpha\p{M;E}$ of $\alpha$-Hölder continuous sections of $E$. In order to define the topology on this space, we recall that we fixed a smooth Riemannian metric on $M$ and choose $\epsilon > 0$ such that for $x,y \in M$ with $d(x,y) \leq 2\epsilon$ there is a unique unit speed geodesic $c_{y,x}$ from $y$ to $x$ in the ball of center $x$ and radius $2 \epsilon$. For such $x$ and $y$, we let $S_{x,y} : E_y \to E_x$ denote the parallel transport \cite[\S 2.B.5]{ghl_geometry} for $\nabla$ along $c_{y,x}$ (which is well-defined since $r >1$). We can then define for $v \in \mathcal{C}^\alpha\p{M;E}$ the best Hölder constant for $v$ as
\begin{equation}\label{eq:best_holder}
\va{v}_{\alpha} \coloneqq \sup_{\substack{x,y \in M \\ d(x,y) \leq \epsilon}} \frac{\va{v(x) - S_{x,y} v(y)}}{d(x,y)^\alpha}.
\end{equation}
A norm on $\mathcal{C}^{\alpha}\p{M;E}$ is then defined by
\begin{equation}
\label{eq:norm_def_holder}
\n{v}_\alpha \coloneqq \va{v}_\infty + \va{v}_{\alpha} \quad \textup{for} \quad v \in \mathcal{C}^{\alpha}\p{M;E}.
\end{equation}
The norm $\n{\cdot}_\alpha$ defines a structure of Banach space on $\mathcal{C}^{\alpha}\p{M;E}$. Other choices of Riemannian metric on $M$ and of covariant derivative on $E$ lead to equivalent norms. Notice that one could also define the best Hölder constant \eqref{eq:best_holder} by using a finite number of trivializations for $E$: this procedure also gives an equivalent norm.

Finally, we will need a slightly less classical space of sections of $E$: the space $\mathcal{C}^{\alpha,k}\p{M;E}$ of sections of $E$ that admit $k$ derivatives in the unstable direction that are $\alpha$-Hölder. Similar classes of regularity have already been introduced in the literature on hyperbolic dynamics, see \cite{delallaveRemarksSobolevRegularity2001}. The point of introducing the spaces $\mathcal{C}^{\alpha,k}\p{M;E}$ is that on the one hand they are large enough so that we can perform the decomposition \eqref{eq:decomposition_X} in these spaces (see \S \ref{section:Y_and_Z}), and on the other hand having derivatives in the unstable direction is essential in Steps 3 and 4 of the S3 algorithm (see \S \ref{section:unstable} and \S \ref{section:convergence}). The class of differentiable functions is in general too small to get the decomposition \eqref{eq:decomposition_X}. Indeed, since the unstable direction is in general only H\"older continuous, asking for a differentiability of $Y$ in \eqref{eq:decomposition_X} would be too restrictive.

We say that a section $v$ of $E$ is differentiable in the unstable direction if its restriction to any unstable manifold is differentiable. In that case, we can define $\nabla_{X} v$, when $X$ is a vector field tangent to $E_u$, by applying the pullback of $\nabla$ onto the unstable manifold to the restriction of $v$. Let us detail this definition. Let $x_0$ be a point in $M$ and $W$ be a local unstable manifold for $\phi$ at $x_0$. Since $W$ is a $\mathcal{C}^r$ manifold and the restriction of $v$ to $W$ is differentiable, we can find a differentiable section $\tilde{v}$ of $E$ that coincides with $v$ on $W$. Then, it makes sense to consider $\nabla_{X} \tilde{v}$ (even if $X$ is not smooth, see \eqref{eq:expression_coordinate}), and since $X$ is tangent to $W$, we see that the value of $\nabla_{X} \tilde{v}$ in the interior of $W$ does not depend on the choice of $\tilde{v}$.  In particular it makes sense to set
\begin{equation*}
\nabla_{X} v (x_0) \coloneqq \nabla_{X} \widetilde{v}(x_0).
\end{equation*}

We can then define the space $\mathcal{C}^{\alpha,k}\p{M;E}$ by induction on $k$, setting $\mathcal{C}^{\alpha,0}\p{M;E} = \mathcal{C}^{\alpha}\p{M;E}$. Then for $k \in \N$, we say that $v \in \mathcal{C}^{\alpha,k+1}\p{M;E}$ if $v \in \mathcal{C}^{\alpha,k}\p{M;E}$, the section $v$ is differentiable in the unstable direction and $\nabla_{X} v \in \mathcal{C}^{\alpha,k}\p{M;E}$, for all $X \in \mathcal{C}^{\alpha,k}\p{M,TM}$ tangent to the unstable direction. If $v \in \mathcal{C}^{\alpha,k}\p{M;E}$, we will sometimes say for short that $v$ is $\mathcal{C}^{\alpha,k}$. We will also write $\mathcal{C}^{\alpha,k}$ instead of $\mathcal{C}^{\alpha,k}\p{M;M \times \C}$. It follows from Remark \ref{remark:different_connection} that, under the assumption of Lemma \ref{lemma:spaces}, the $\mathcal{C}^{\alpha,k}$ regularity is preserved by basic algebraic operations (such as product of endomorphisms or scalar product). Moreover, the formulae given in Remark \ref{remark:different_connection} remain valid when differentiating a $\mathcal{C}^{\alpha,k}$ sections in the unstable direction, as can be seen by replacing $\mathcal{C}^{\alpha,k}$ sections by smooth sections that coincide with them on a local unstable manifold.

We will say that the unstable direction $E_u$ is $\mathcal{C}^{k,\alpha}$ if the orthogonal projector on $E_u$ is $\mathcal{C}^{\alpha,k}$ as a section of the bundle of endomorphism of $TM$, that is if it belongs to $\mathcal{C}^{\alpha,k}\p{M,\End(TM)}$. Using Gram--Schmidt process, one can see that $E_u$ is $\mathcal{C}^{\alpha,k}$ if and only if for every $x_0 \in M$, there are $\mathcal{C}^{\alpha,k}$ vector fields $X_1,\dots,X_{d_u}$ such that $X_1(x),\dots,X_{d_u}(x)$ is a basis of $E_u(x)$ for every $x$ near $x_0$ (in particular this definition does not depend on the choice of Riemannian metric on $M$). By a compactness argument, we also see that $E_u$ is $\mathcal{C}^{\alpha,k}$ if and only if there is a family $X_1,\dots,X_D$ of $\mathcal{C}^{\alpha,k}$ vector fields tangent to $E_u$ such that, for every $x \in M$, the vectors $X_1(x),\dots,X_D(x)$ span $E_u(x)$ (notice that in general we expect $D > d_u$). A consequence of this fact is that if $k+ \alpha < r$ and $E_u$ is $\mathcal{C}^{\alpha,k-1}$ then the restriction to a local unstable manifold of a $\mathcal{C}^{\alpha,k}$ function is $\mathcal{C}^{k + \alpha}$.

\begin{remark}
We gave a definition of the space $\mathcal{C}^{\alpha,k}\p{M;E}$ for any $k \in \N$ and $0 \leq \alpha < 1$. However, when $k \geq 1$, it is not clear \emph{a priori} that this definition does not depend on the choice of covariant derivative on $E$. Since the unstable direction $E_u$ is not smooth, it is not even clear that the space $\mathcal{C}^{\alpha,k}\p{M;E}$ is non-trivial. It follows from Lemma \ref{lemma:spaces} below that the space $\mathcal{C}^{\alpha,k}\p{M;E}$ is well-behaved as soon as $k + \alpha \leq r$ and $E_u$ is $\mathcal{C}^{\alpha,k-1}\p{M;E}$. It might be possible to deduce from standard techniques in hyperbolic dynamical systems that this condition is satisfied if $k < r-1$ and $\alpha$ is small enough (see for instance \cite[Theorem 2]{delallaveRemarksSobolevRegularity2001} and references therein). However, we will prove this regularity of $E_u$ by a bootstrap argument using Proposition \ref{prop:bound_spectral_radius}, see Lemma \ref{lemma:regularity}.
\end{remark}

If $k \geq 1$ and $E_u$ is $\mathcal{C}^{\alpha,k-1}$, we define a norm on $\mathcal{C}^{\alpha,k}(M,E)$ in the following way. Fix a family $X_1,\dots,X_D$ of $\mathcal{C}^{\alpha,k-1}$ vector fields tangent to $E_u$ such that, for every $x \in M$, the vectors $X_1(x),\dots,X_D(x)$ span $E_u(x)$. We endow $\mathcal{C}^{\alpha,k}\p{M;E}$ with the norm
\begin{equation}\label{eq:def_norm}
\n{v}_{\alpha,k} = \sum_{j = 0}^{k} \sum_{1 \leq i_1, \dots, i_j \leq D} \|\nabla_{X_{i_1}} \dots \nabla_{X_{i_j}} v\|_{\alpha}
\end{equation}

Let us mention here some elementary properties of the spaces $\mathcal{C}^{\alpha,k}$ that will be useful later.

\begin{lemma}\label{lemma:spaces}
Let $k \in \N$ and $0 \leq \alpha < 1$ be such that $k+\alpha \leq r$. If $k \geq 1$, assume in addition that $E_u$ is $\mathcal{C}^{\alpha,k-1}$. Then the following properties hold:
\begin{enumerate}[label=(\roman*)]
\item if $v$ is a $\mathcal{C}^{k + \alpha}$ section of $E$, then $v \in \mathcal{C}^{\alpha,k}\p{M;E}$;
\item the space $\mathcal{C}^{\alpha,k}\p{M;E}$ does not depend on the choice of the covariant derivative $\nabla$ (up to taking an equivalent norm);
\item $\mathcal{C}^{\alpha,k}\p{M;E}$ is a Banach space;
\item if $k \geq 1$, the vector fields $X_1,\dots,X_D$ are $\mathcal{C}^{\alpha,k-1}$, tangent to $E_u$ and span it, and $Y \in \mathcal{C}^{\alpha,k-1}\p{M;TM}$ is parallel to the unstable direction, then there are $a_1,\dots,a_D \in \mathcal{C}^{\alpha,k-1}$ such that $Y = \sum_{j = 1}^D a_j X_j$. Moreover, there is a constant $C > 0$ that does not depend on $Y$ such that $\va{a_j}_\infty \leq C \va{Y}_\infty$ for $j = 1,\dots,D$.
\end{enumerate}
\end{lemma}

\begin{proof}
\noindent (i) The restriction of $v$ to any unstable manifold is $\mathcal{C}^k$. If $0 \leq \ell \leq k$ and $X_1,\dots,X_\ell$ are $\mathcal{C}^{\alpha,k-1}$ vector fields tangent to $E_u$, we can consequently define $\nabla_{X_\ell} \dots \nabla_{X_1} v$. Using the expression \eqref{eq:expression_coordinate} for $\nabla$ in local coordinates, we see that $\nabla_{X_\ell} \dots \nabla_{X_1} v$  may be expressed in these coordinates as a polynomial in the derivatives of $v$ of order at most $\ell$ and the $\nabla_{X_{i_m}} \dots \nabla_{X_{i_1}} X_{i_0}$'s for $0 \leq m \leq \ell - 1$ and $0 \leq i_0 \leq \dots \leq i_m \leq \ell$, with $\mathcal{C}^{r - \ell}$ coefficients (to actually prove this formula when $\ell \geq 1$, consider the restriction of $X_1,\dots,X_\ell$ to a local unstable manifold and extend them to $\mathcal{C}^{k-1}$ vector fields near that curve). Hence, we see, by induction on $\ell$, that $\nabla_{X_\ell} \dots \nabla_{X_1} v$ is indeed $\alpha$-Hölder.

\noindent (ii) The proof is an induction on $k$, the case $k = 0$ being well-known. Assume that $k \geq 1$ and that the result holds for $k-1$. Let $\widetilde{\nabla}$ be another $\mathcal{C}^r$ covariant derivative on $E$. Then, it follows from \eqref{eq:expression_coordinate} that there is a $\mathcal{C}^{r-1}$ section $\gamma$ of the vector bundle $T^*M \otimes E^* \otimes E$ such that for $X$ a vector field and $v$ a section of $E$, we have
\begin{equation*}
\widetilde{\nabla}_X v = \nabla_{X} v + \gamma\p{X,v}.
\end{equation*}
If $v$ belongs to the space $\mathcal{C}^{\alpha,k}\p{M;E}$ defined using the covariant derivative $\nabla$ and $X$ is $\mathcal{C}^{\alpha,k-1}$ and tangent to $E_u$, we know that $\nabla_X v$ and $v$ belong to $\mathcal{C}^{\alpha,k-1}\p{M;E}$, whose definition does not depend on $\nabla$. Applying Leibniz rule \eqref{eq:Leibniz} and point (i) to $\gamma$, we see that $\widetilde{\nabla}_{X} v$ belongs to $\mathcal{C}^{\alpha,k-1}\p{M;E}$. Consequently, $v$ also belongs to the space $\mathcal{C}^{\alpha,k}\p{M;E}$ defined using the covariant derivative $\widetilde{\nabla}$. It follows that the space $\mathcal{C}^{\alpha,k}\p{M;E}$ does not depend on the choice of $\nabla$.

\noindent (iii) Let $(v_n)_{n \in \N}$ be a Cauchy sequence in $\mathcal{C}^{\alpha,k}\p{M;E}$. Let $X_1,\dots,X_D$ be as in the definition \eqref{eq:def_norm} of the norm on $\mathcal{C}^{\alpha,k}\p{M;E}$. Since $\mathcal{C}^{\alpha}\p{M;E}$ is a Banach space (with norm \eqref{eq:norm_def_holder}), for $0 \leq  j \leq k$ and $1 \leq i_1,\dots,i_k \leq D$, the sequence $(\nabla_{X_{i_1}} \dots \nabla_{X_{i_j}} v_n)_{n \in \N}$ converges to an element $u_{i_1,\dots,i_j}$ of $\mathcal{C}^{\alpha}\p{M;E}$. Hence, we only need to prove that the restriction of $u_0$ to any unstable submanifold is $\mathcal{C}^k$ and that $\nabla_{X_{i_1}} \dots \nabla_{X_{i_j}} u_0 = u_{i_1,\dots,i_j}$. This follows by restricting to local unstable manifolds and applying the fact that the space of $\mathcal{C}^{k}$ sections of a $\mathcal{C}^k$ vector bundle over a manifold with boundary is a Banach space. 

\noindent (iv) By compactness, one can work locally, near a point $x_0$. Then, up to reordering $X_1,\dots,X_D$, one can assume that $X_1(x),\dots,X_{d_u}(x)$ is a basis of $E_u(x)$ for $x$ near $x_0$ (recall that \emph{a priori} $D > d_u$). Then, for such an $x$, we can write $Y(x) = \sum_{j = 1}^{d_u} a_j(x) X_j(x)$, where the $a_j(x)$ are obtained by solving the system of equations $\langle Y(x), X_j(x) \rangle = \sum_{i = 1}^{d_u} a_i(x) \langle X_i(x),X_j(x) \rangle, 1\leq j\leq D$. Since this system may be solved using Cramer's rule (the determinant of the associated matrix is non-zero), we obtain that the $a_j$'s are $\mathcal{C}^{\alpha,k-1}$ with their suprema controlled by the suprema of $Y$.
\end{proof}

\begin{remark}\label{remark:coefficients}
Notice that the method used to construct the coefficients $a_j$'s in point (iv) of Lemma \ref{lemma:spaces} allow to compute explicitly the $a_j$'s from the knowledge of $Y$ and $X_1,\dots,X_D$. Moreover, if $X$ is a vector field parallel to the unstable direction and $k \geq 2$, one may compute the $Xa_j$'s if they know $\nabla_X Y,\nabla_X X_1,\dots,\nabla_X X_D$. Indeed, using the Levi--Civita connection we see that the $X a_j$'s satisfy the system of equations (using the same notation as in the proof of Lemma \ref{lemma:spaces})
\begin{equation*}
\begin{split}
\sum_{i = 1}^{d_u} X a_i \langle X_i, X_j \rangle & = \langle \nabla_X Y, X_j \rangle + \langle Y,\nabla_X X_j \rangle \\ & \qquad \qquad - \sum_{i = 1}^{d_u} a_i \langle \nabla_X X_i, X_j \rangle - \sum_{i = 1}^{d_u} a_i \langle X_i, \nabla_X X_j \rangle.
\end{split}
\end{equation*}
This system can be solved using Cramer's rule. Notice also that the $a_j$ depends linearly on $Y$ and that, iterating the procedure above, we find that $\n{a_j}_{\alpha,k-1} \leq C \n{Y}_{\alpha,k-1}$ for some constant $C$ that does not depend on $Y$. Consequently, an error in the determination of $Y$ and its derivatives in the unstable direction will produce an error of the same order of magnitude in the determination of the $a_j$'s and their derivatives in the unstable direction.
\end{remark}

In order to bound the spectral radius of the transfer operator $\mathcal{L}_{\Phi}$ defined by \eqref{eq:general_transfer_operator} acting on $\mathcal{C}^{\alpha,k}\p{M;E}$, let us introduce the quantity
\begin{equation*}
\theta = \theta(\Phi) \coloneqq \limsup_{n \to + \infty} \p{\sup_{x \in M} \sup_{v \in E_x \setminus \set{0}} \frac{\va{\Phi^n v}}{\va{v}}}^{\frac{1}{n}}.
\end{equation*}
We will also need the maximal expansion rate for $\phi^{-1}$:
\begin{equation}\label{eq:def_Lambda}
\Lambda \coloneqq \limsup_{n \to + \infty} \p{\sup_{x \in M} \sup_{v \in T_x M \setminus \set{0}} \frac{\va{\mathrm{d}_x \phi^{-n} v}}{\va{v}}}^{\frac{1}{n}}.
\end{equation}

\begin{remark}\label{remark:regularity_lift}
Let us notice that a lift $\Phi$ of $\phi$ identifies with a section of the bundle $\phi^*E \otimes E^*$ or $E \otimes (\phi^{-1})^* E^*$. Consequently, it makes sense to say that $\Phi$ is $\mathcal{C}^{\alpha,k}$.
\end{remark}

We can now state the main result of this section.

\begin{prop}\label{prop:bound_spectral_radius}
Let $0 \leq \alpha < 1$ and $k \leq r$. Assume that $\Phi$ is $\mathcal{C}^{\alpha,k}$ (in the sense of Remark \ref{remark:regularity_lift}), that $k + \alpha \leq r$  and that $E_u$ is $\mathcal{C}^{\alpha,\max(k-1,0)}$. Then $\mathcal{L}_{\Phi}$ induces a bounded operator on $\mathcal{C}^{\alpha,k}\p{M;E}$ with spectral radius less than $\Lambda^\alpha \theta$.
\end{prop}

Before proving Proposition \ref{prop:bound_spectral_radius}, let us discuss how we will use it in the context of the S3 algorithm.

\begin{remark}\label{remark:solving_equations}
In order to compute $Y,V$ and $\rho_Y$ in Steps 2 and 3 of the algorithm, we will need to solve for $v$ equations of the form
\begin{equation}\label{eq:equation}
(I - \mathcal{L}_{\Phi}) v= w,
\end{equation}
where $w$ is a section of a vector bundle $E$ over $M$ and $\Phi : E \to E$ a lift of $\phi$, see Remarks \ref{remark:step_2} and \ref{remark:compute_rho}. If $\Phi$ and $w$ are continuous and $\theta(\Phi) <1$, then it follows from Proposition \ref{prop:bound_spectral_radius} that \eqref{eq:equation} has a unique continuous solution $v$ which is given by
\begin{equation}\label{eq:solution_equation}
v = (I - \mathcal{L}_{\Phi})^{-1} w = \sum_{n \geq 0} \mathcal{L}_{\Phi}^n w.
\end{equation}
It will be useful to notice that if $\Phi$ and $w$ are Hölder-continuous, then so is $v$, since we can always find a small $\alpha > 0$ such that $\theta(\Phi) \Lambda^\alpha < 1$. Similarly, if $\Phi$ and $w$ are $\mathcal{C}^{\alpha,k}$ then so is $v$.

When actually implementing the S3 algorithm, equations of the form \eqref{eq:equation} are solved approximately by choosing a large $L$ and using as a solution
\begin{equation*}
v_L = \sum_{n = 0}^{L - 1} \mathcal{L}_{\Phi}^n w \simeq v.
\end{equation*}
Since the series in the right hand side of \eqref{eq:solution_equation} converges exponentially fast, we see that the error we make when approximating $v$ by $v_L$ is exponentially small in $L$. Moreover, we see from the definition \eqref{eq:general_transfer_operator} of the operator $\mathcal{L}_{\Phi}$ that, in order to compute the value of $v_L$ at a point $x \in M$, one only needs to know the value of $w$ at $x, \phi^{-1} x, \dots, \phi^{-L + 1} x$, which makes the computation doable.
\end{remark}

\begin{remark}\label{remark:solving_derivative}
When the solution $v$ to \eqref{eq:equation} is $\mathcal{C}^{\alpha,1}$, we will sometimes need to compute derivatives of the form $\nabla_{X} v$, with $X$ parallel to the unstable direction (see Remark \ref{remark:compute_rho}). To do so, one can work as in the proof of Proposition \ref{prop:bound_spectral_radius} in the case $k \geq 1$ below. Let $X_1,\dots,X_D$ be $\mathcal{C}^{\alpha,1}$ vector fields on $M$, parallel to the unstable direction and such that $X_1(x),\dots,X_D(x)$ spans $E_u(x)$ for every $x \in M$. Then, let $\Theta$ be the map from $\mathcal{C}^{\alpha,1}\p{M;E}$ to $\mathcal{C}^{\alpha,0}\p{M; E^{\oplus (D+1)}}$ given by $\Theta(u) = (u \oplus \nabla_{X_1} u \oplus \dots \oplus \nabla_{X_D} u)$. In the proof of Proposition \ref{prop:bound_spectral_radius} in the case $k \geq 1$ below, we find an integer $m \geq 1$ and a lift $\hat{\Phi}_m$ of $\phi^m$ such that $\Theta \circ \mathcal{L}_{\Phi}^m = \mathcal{L}_{\hat{\Phi}_m} \circ \Theta$ and the spectral radius of $\mathcal{L}_{\hat{\Phi}_m}$ on $\mathcal{C}^{\alpha,0}\p{M; E^{\oplus (D+1)}}$ is less than $1$ (provided $\theta(\Phi) < 1$).

We have consequently
\begin{equation*}
(I - \mathcal{L}_{\hat{\Phi}_m}) \Theta(v) = \Theta((I - \mathcal{L}_{\Phi}^m)v) = \Theta((I + \mathcal{L}_{\Phi} + \dots + \mathcal{L}_{\Phi}^{m-1})w),
\end{equation*}
so that
\begin{equation*}
\begin{split}
\Theta(v) & = (I - \mathcal{L}_{\hat{\Phi}_m})^{-1} \Theta((I + \mathcal{L}_{\Phi} + \dots + \mathcal{L}_{\Phi}^{m-1})w) \\ & = \sum_{n \geq 0} \mathcal{L}_{\hat{\Phi}_m}^n \Theta((I + \mathcal{L}_{\Phi} + \dots + \mathcal{L}_{\Phi}^{m-1})w).
\end{split}
\end{equation*}
We can use this formula to find an approximate value of $\Theta(v)$ as we did in Remark \ref{remark:solving_equations} for $v$. From the knowledge of $\nabla_{X_1} v,\dots,\nabla_{X_D} v$ we can then deduce the value of $\nabla_X v$ for any vector field $X$ parallel to the unstable direction.
\end{remark}

\begin{example}\label{example:basic}
Let us illustrate the method from Remarks \ref{remark:solving_equations} and \ref{remark:solving_derivative} through a very basic example. Assume that $\phi$ is a CAT map of the two-dimensional torus (a linear example of Anosov diffeomorphism) and that we try to solve for $v$ the equation
\begin{equation}\label{eq:example}
v - \frac{1}{2} v \circ \phi^{-1} = w,
\end{equation}
where $w$ is a smooth function. This is just the equation \eqref{eq:equation} in the context of Example \ref{example:scalar} with $h$ identically equal to $1/2$. By a von Neumann series argument, we find a continuous solution $v$ to \eqref{eq:example} given by
\begin{equation*}
v = \sum_{n \geq 0} 2^{-n} w \circ \phi^{-n}.
\end{equation*}
We want now to investigate the regularity of $v$. Let $\lambda_u > 1$ and $\lambda_s = \lambda_u^{-1}$ be the eigenvalues of $\phi$. Differentiating in the stable direction, we see that each term is multiplied by a factor $\lambda_u^n$, so that $v$ will not be differentiable (unless $\lambda_u$ is less than $2$). This can be mitigated by choosing $\alpha > 0$ such that $\lambda_u^\alpha < 2$ and noticing that $v$ is $\alpha$-Hölder. Differentiating in the unstable direction, the situation is much better since each term is multiplied by a factor $\lambda_s^n$, so that we can see that $v$ is indeed differentiable in the unstable direction. The proof of Proposition \ref{prop:bound_spectral_radius} is an adaptation of this observation.
\end{example} 

\begin{remark}
A natural question regarding the implementation of the S3 algorithm is whether the inevitable computational error will affect the quality of the approximation of the solution to an equation of the form \eqref{eq:equation}. A small error when computing the action of $\Phi$ on the fiber amounts to replacing the operator $\mathcal{L}_{\Phi}$ by an operator which is close in the norm topology of operators from $\mathcal{C}^{\alpha,k}\p{M;E}$ to itself. Thus, the resulting error is still small at the level of the resolvent $(I - \mathcal{L}_\Phi)^{-1}$.

An error when computing the action of $\phi$ on a point $x$ is more serious. Indeed, such a perturbation of the dynamics induces a perturbation of the operator $\mathcal{L}_{\Phi}$ which may not be small in the norm topology of operators from $\mathcal{C}^{\alpha,k}\p{M;E}$ to itself (the perturbation is small in the strong operator topology but this is not enough to get spectral stability). However, if we see $\mathcal{L}_{\Phi}$ as an operator from $\mathcal{C}^{\alpha,k}$ to $\mathcal{C}^{\beta,k}$ with $\beta < \alpha$, then the perturbation becomes small. One can then use the machinery developed by Gouëzel, Keller and Liverani in the context of linear response to get some kind of spectral stability \cite{Ke-Li, GLK} (see \cite[A.3]{Bal2} for a textbook presentation). In particular, since we are only working away from the spectrum, the resolvent of $\mathcal{L}_{\Phi}$ is well approximated by the resolvent of its perturbation, and so we can expect that the computational errors do not affect the result of the S3 algorithm too much.
\end{remark}

We will prove Proposition \ref{prop:bound_spectral_radius} first in the case $k = 0$.

\begin{proof}[Proof of Proposition \ref{prop:bound_spectral_radius} in the case $k = 0$]
The case $\alpha = 0$ is elementary. It is also a consequence of the much more general statement \cite[Theorem 8.15]{chiconeEvolutionSemigroupsDynamical1999}. Let us consequently assume that $0 < \alpha < 1$.

Let $\widetilde{\Lambda} > \Lambda$ and $\tilde{\theta} > \theta$. Choose some number $\delta \ll 1$, and let $v \in \mathcal{C}^{\alpha,0}\p{M;E}$ and $n \in \N$. From the definition of $\theta$, we see that there is a constant $C > 0$, that does not depend on $n$, such that, for every $x \in M$, the operator norm of $\Phi^n$ from $E_x$ to $E_{\phi^n x}$ is less than $C \tilde{\theta}^n$. It follows that 
\begin{equation}\label{eq:estimate_sup_norm}
\va{\mathcal{L}_{\Phi}^n v}_{\infty} \leq C \tilde{\theta}^n \va{v}_\infty.
\end{equation}
We want to estimate the semi-norm $\va{\mathcal{L}_{\Phi}^n v}_{\alpha}$. Let $x,y \in M$ be such that $d(x,y) \leq \epsilon$. If $d(x,y) \geq \widetilde{\Lambda}^{-n} \delta$, then we have, for some $C > 0$ that may vary from one line to another (but does not depend on $n$ nor $v$),
\begin{equation*}
\begin{split}
\va{\Phi^n v(\phi^{-n} x) - S_{x,y}\Phi^{n} v(\phi^{-n} y)} & \leq C \p{\va{\Phi^n v(\phi^{-n} x)} + \va{\Phi^n v(\phi^{-n} y)}} \\
   & \leq C \tilde{\theta}^n \va{v}_\infty \leq C \delta^{- \alpha} (\tilde{\theta} \widetilde{\Lambda}^\alpha)^n d(x,y)^\alpha.
\end{split}
\end{equation*}
Let us now consider the case $d(x,y) \leq \widetilde{\Lambda}^{-n} \delta$. Notice that there is a constant $C > 0$ such that the derivative of $\phi^{-n}$ is bounded by $C \widetilde{\Lambda}^{n}$. Hence we have
\begin{equation*}
d(\phi^{-n} x,\phi^{-n} y) \leq C \widetilde{\Lambda}^n d(x,y) \leq C \delta.
\end{equation*}
In particular, we have $d(\phi^{-n} x,\phi^{-n} y) \leq \epsilon$, provided $\delta$ is small enough, and we can estimate
\begin{equation*}
\begin{split}
& \va{\Phi^n v(\phi^{-n} x) - S_{x,y}\Phi^{n} v(\phi^{-n} y)} \\  & \qquad \leq \va{\Phi^n \p{v(\phi^{-n}x) - S_{\phi^{-n}x, \phi^{-n} y}v(\phi^{-n} y)}}  \\ & \qquad \qquad \qquad \qquad \qquad + \va{\Phi^n S_{\phi^{-n}x, \phi^{-n} y}v(\phi^{-n} y) - S_{x,y} \Phi^n v(\phi^{-n} y)} \\
   & \qquad \leq C \tilde{\theta}^n \va{v}_{\alpha} d(\phi^{-n} x, \phi^{-n} y)^\alpha \\ & \qquad \qquad \qquad \qquad \qquad + \va{\Phi^n S_{\phi^{-n}x, \phi^{-n} y}v(\phi^{-n} y) - S_{x,y} \Phi^n v(\phi^{-n} y)} \\
   & \qquad \leq C \va{v}_\alpha (\tilde{\theta} \widetilde{\Lambda}^\alpha)^n d(x,y)^\alpha + \va{\Phi^n S_{\phi^{-n} x, \phi^{-n} y}v(\phi^{-n} y) - S_{x,y} \Phi^n v(\phi^{-n} y)}. 
\end{split}
\end{equation*}
We want then to estimate the operator norm of $\Phi^n S_{\phi^{-n }x, \phi^{-n} y} - S_{x,y} \Phi^n$, which is a linear application from $E_{\phi^{-n}y}$ to $E_x$. To do so, we write
\begin{equation*}
\begin{split}
& \Phi^n S_{\phi^{-n}x, \phi^{-n} y} - S_{x,y} \Phi^n \\ & \qquad \qquad = \sum_{\ell = 0}^{n-1} \Phi^\ell S_{\phi^{- \ell} x, \phi^{- \ell} y}\p{S_{\phi^{- \ell} y, \phi^{- \ell} x} \Phi S_{\phi^{- \ell - 1} x, \phi^{- \ell - 1} y} - \Phi} \Phi^{n - \ell - 1}.
\end{split}
\end{equation*}
In each term, we may bound the operator norm of $\Phi^\ell$ by $C \tilde{\theta}^\ell$, the operator norm of $\Phi^{n - \ell - 1}$ by $C \tilde{\theta}^{n- \ell}$ and the operator norm of $S_{\phi^{- \ell} x, \phi^{- \ell} y}$ by a constant. Notice then that, for $\ell = 0,\dots ,n-1$, the distance between $\phi^{- \ell - 1} x$ and $\phi^{- \ell - 1} y$ is less than $C \widetilde{\Lambda}^\ell d(x,y)$, which is less than $\epsilon$, provided $\delta$ is small enough. Consequently, it follows from our hypothesis on $\Phi$ that the operator norm of 
\begin{equation*}
S_{\phi^{- \ell} y, \phi^{- \ell} x} \Phi S_{\phi^{- \ell - 1} x, \phi^{- \ell - 1} y} - \Phi
\end{equation*}
is less than $C \widetilde{\Lambda}^{\alpha \ell} d(x,y)^\alpha$. The operator norm of $\Phi^n S_{\phi^{-n }x, \phi^{-n} y} - S_{x,y} \Phi^n$ is consequently less than $ C (\tilde{\theta} \widetilde{\Lambda}^\alpha)^n d(x,y)^\alpha $ with $C$ that does not depend on $n$. 

Hence, we find that for some $C > 0$ that does not depend on $n$ nor $v$, we have
\begin{equation*}
\begin{split}
\va{\mathcal{L}_\Phi^n v}_\alpha \leq C (\tilde{\theta} \widetilde{\Lambda}^\alpha)^n \n{v}_\alpha.
\end{split}
\end{equation*}
Recalling \eqref{eq:estimate_sup_norm}, it follows that the spectral radius of $\mathcal{L}_\Phi$ acting on $\mathcal{C}^{\alpha}\p{M;E}$ is less than $\tilde{\theta} \widetilde{\Lambda}^\alpha$. Since $\tilde{\theta} > \theta$ and $\widetilde{\Lambda} > \Lambda$ are arbitrary, the result follows.
\end{proof}

In order to deduce the general case in Proposition \ref{prop:bound_spectral_radius} from the case $k = 0$, we will need the following Lemma, which follows from the expression \eqref{eq:expression_coordinate} for $\nabla$ in local coordinates and the chain rule.

\begin{lemma}\label{lemma:chain_rule}
Let $0 \leq \alpha < 1$ and $1 \leq k \leq r$. Assume that $\Phi$ is $\mathcal{C}^{\alpha,k}$, that $k + \alpha \leq r$  and that $E_u$ is $\mathcal{C}^{\alpha,k-1}$. Then, there is a $\mathcal{C}^{\alpha,k-1}$ section $\gamma$ of the bundle $T^* M \otimes (\phi^{-1})^* E^* \otimes E$ such that, for every $v \in \mathcal{C}^{\alpha,k}\p{M;E}$ and $X \in \mathcal{C}^{\alpha,k-1}\p{M;TM}$ parallel to the unstable direction, we have
\begin{equation}\label{eq:chain_rule}
\nabla_X \p{\mathcal{L}_{\Phi} v} = \mathcal{L}_{\Phi}\p{\nabla_{\phi_*^{-1} X} v} + \gamma(X,v \circ \phi^{-1}).
\end{equation}
\end{lemma}

\begin{proof}
As in Remark \eqref{remark:regularity_lift}, let us identify the lift $\Phi$ with a section of the bundle $E \otimes (\phi^{-1})^* E^*$. We define the covariant derivative on this vector bundle as in Remark \ref{remark:different_connection}. With this identification and using \eqref{eq:Leibniz_endomorphism}, we have
\begin{equation*}
\begin{split}
\nabla_X \p{\mathcal{L}_{\Phi} v} & = \nabla_X \p{\Phi (v \circ \phi^{-1})}  = \p{\nabla_X \Phi}(v \circ \phi^{-1}) + \Phi\p{\nabla_X (v \circ \phi^{-1})} \\
   & = \Phi\p{\p{\nabla_{\phi_*^{-1} X} v} \circ \phi^{-1}} + \p{\nabla_X \Phi}(v \circ \phi^{-1}) \\
   & = \mathcal{L}_\Phi\p{\nabla_{\phi_*^{-1} X} v} +\p{\nabla_X \Phi}(v \circ \phi^{-1}).
\end{split}
\end{equation*}
Notice that on the second line we used the particular choice of covariant derivative on $(\phi^{-1})^* E$ from Remark \ref{remark:different_connection}. It follows from \eqref{eq:expression_coordinate} that the term $\p{\nabla_X \Phi}(v \circ \phi^{-1})$ is of the form $\gamma(X,v \circ \phi^{-1})$ as expected. This proof can also be written entirely in coordinates using \eqref{eq:expression_coordinate}.
\end{proof}

We can now use Lemma \ref{lemma:chain_rule} to deduce the general case in Proposition \ref{prop:bound_spectral_radius} from the case $k=0$.

\begin{proof}[Proof of Proposition \ref{prop:bound_spectral_radius} in full generality]
The plan is to reduce to the case $k = 0$. Let $X_1,\dots,X_D$ be as in the definition \eqref{eq:def_norm} of the norm on $\mathcal{C}^{\alpha,k}\p{M;E}$. Let us consider a new vector bundle over $M$:
\begin{equation*}
F = \bigoplus_{j = 0}^{k} F_j, \textup{ where } F_j = \bigoplus_{1 \leq i_1,\dots,i_j \leq D} E,
\end{equation*}
with the convention that $F_0 = E$. We define then the map
\begin{equation*}
\begin{array}{ccccc}
\Theta & : & \mathcal{C}^{\alpha,k}\p{M;E} & \to & \mathcal{C}^{\alpha}\p{M;F} \\
 & & v & \mapsto & \p{\Theta_j(v)}_{0 \leq j \leq k}.
\end{array}
\end{equation*}
where
\begin{equation*}
\Theta_j(v) = \p{\nabla_{X_{i_1}} \dots \nabla_{X_{i_j}} v}_{ 1 \leq i_1,\dots,i_j \leq D}.
\end{equation*}
It follows from Lemma \ref{lemma:spaces} that the map $\Theta$ induces a continuous linear isomorphism between $\mathcal{C}^{\alpha,k}\p{M;E}$ and its (closed) image in $\mathcal{C}^{\alpha}\p{M;F}$. Let then $m \geq 1$ be a large integer, for $i = 1,\dots,D$, let $a_{i1}^m,\dots,a_{iD}^m$ be $\mathcal{C}^{\alpha,k-1}$ functions given by Lemma \ref{lemma:spaces} such that
\begin{equation}\label{eq:decompostion_Xi}
\phi_*^{-m} X_i = \sum_{j = 1}^D a_{ij}^m X_j.
\end{equation}
Since $X_i$ is parallel to the unstable direction, we see that, for $i,j = 1,\dots,D$, we have 
\begin{equation}\label{eq:borne_aij}
|a_{ij}^m| \leq C \delta^m,
\end{equation}
where $C > 0$ and $0 < \delta <1$ do not depend on $m$. Hence, applying Lemma \ref{lemma:chain_rule} iteratively, we see that for $j = 1,\dots,k$ and $1 \leq i_1,\dots, i_j \leq D$, we have
\begin{equation}\label{eq:calcul_derivative}
\begin{split}
& \nabla_{X_{i_1}} \dots \nabla_{X_{i_j}} \p{\mathcal{L}_{\Phi}^{m} v} \\ & \qquad = \sum_{1 \leq p_1,\dots ,p_j \leq D} a_{i_1 p_1}^m \circ \phi^{-m} \dots a_{i_j p_j}^m \circ \phi^{-m}  \mathcal{L}_{\Phi}^m\p{\nabla_{X_{p_1}} \dots \nabla_{X_{p_j}}v} \\& \qquad \qquad \qquad \quad + \Gamma_{i_1,\dots,i_j}^{m}(v \circ \phi^{-m},\Theta_1(v) \circ \phi^{-m},\dots,\Theta_{j-1}(v) \circ \phi^{-m}),
\end{split}
\end{equation}
where $\Gamma_{i_1,\dots,i_j}^{m}$ is an $\alpha$-Hölder section of the vector bundle $(\phi^{-m})^{*} F_0^* \otimes \dots \otimes (\phi^{-m})^* F_{j-1}^* \otimes E$. Now, let $A_{j,m}$ be the section of $\End(F_j)$ whose action is given by the matrix $(a_{i_1 p_1}^m \circ \phi^{-m} \dots a_{i_j p_j}^m \circ \phi^{-m})_{\substack{1 \leq i_1,\dots,i_j \leq D \\ 1 \leq p_1,\dots,p_j \leq D }}$, with the convention that $A_{0,m}$ is the identity, and $\mathcal{A}_m$ be the section of $\End(F)$ that acts diagonally \emph{via} the $A_{j,m}$'s: $\mathcal{A}_m(u_j)_{0 \leq j \leq k} = (A_{j,m} u_j)_{0 \leq j \leq k}$. Let $\widetilde{\Phi}^{m} : F \to F$ be the lift of $\phi^m$ that acts like $\Phi^m$ in each component. According to \eqref{eq:calcul_derivative}, there is another lift $R_m: F \to F$ of $\phi^m$ which is striclty lower triangular with respect to the decomposition $F = \bigoplus_{j = 0}^k F_j$ such that if we set
\begin{equation*}
\hat{\Phi}_m = \mathcal{A}_m \widetilde{\Phi}^m  + R_m,
\end{equation*}
then $\mathcal{L}_{\hat{\Phi}_m} \circ \Theta = \Theta \circ \mathcal{L}_{\Phi^m} = \Theta \circ \mathcal{L}_{\Phi}^m$. The lower triangular term $R_m$ is here to take into account $\Gamma_{i_1,\dots,i_j}^{m}$ in \eqref{eq:calcul_derivative}. We want now to apply the case $k = 0$ in Proposition \ref{prop:bound_spectral_radius} to the operator $\mathcal{L}_{\hat{\Phi}_m}$. We start by noticing that this lift of $\phi^m$ is $\alpha$-Hölder. Then, replacing $\phi$ by $\phi^m$ in \eqref{eq:def_Lambda}, we replace $\Lambda$ by at most $\Lambda^m$. Hence, the spectral radius of $\mathcal{L}_{\hat{\Phi}_m}$ acting on $\mathcal{C}^{\alpha}(M;F)$ is at most $\theta(\hat{\Phi}_m) \Lambda^{\alpha m}$. From $\mathcal{L}_{\hat{\Phi}_m} \circ \Theta = \Theta \circ \mathcal{L}_{\Phi}^m$, we deduce that the spectral radius of $\mathcal{L}_{\Phi}$ acting on $\mathcal{C}^{\alpha,k}\p{M;E}$ is at most
\begin{equation*}
\theta(\hat{\Phi}_m)^{\frac{1}{m}} \Lambda^{\alpha}.
\end{equation*}
Since, with respect to the decomposition $F = \bigoplus_{j = 0}^k F_j$, the operator $R_m$ is strictly lower triangular and $\mathcal{A}_m \widetilde{\Phi}^m$ is diagonal, we see that $\theta(\hat{\Phi}_m) = \theta(\mathcal{A}_m \widetilde{\Phi}^m)$. From the bound on the $a_{ij}$'s in \eqref{eq:decompostion_Xi}, we see that the operator norm of $\mathcal{A}_m$ (in each fiber) is bounded uniformly in $m$. Since $\widetilde{\Phi}^m$ acts diagonally like $\Phi^m$, it follows that 
\begin{equation*}
\limsup_{m \to + \infty}\theta(\hat{\Phi}_m)^{\frac{1}{m}} \leq \theta(\Phi),
\end{equation*}
which ends the proof of the proposition.
\end{proof}

Before going into the details of the S3 algorithm, let us give a first application of Proposition \ref{prop:bound_spectral_radius}. 

\begin{lemma}\label{lemma:regularity}
Let $k$ be an integer strictly less than $r-1$. Then there is $0 < \alpha < 1$ such that $E_u$ is $\mathcal{C}^{\alpha,k}$.
\end{lemma}

The proof is based on a higher dimensional version of the argument from \cite[Lemma 8.7]{chandramoorthyEfficientComputationLinear2021}. While this result is not surprising to the specialists in hyperbolic dynamics, the advantage of our proof is that it relies on an application of Proposition \ref{prop:bound_spectral_radius} and thus makes possible to apply the strategy from Remark \ref{remark:solving_equations} to compute explicitly ``the derivative of $E^u$ in the unstable direction'' (see Remark \ref{remark:unstable_basis}).

\begin{proof}[Proof of Lemma \ref{lemma:regularity}]
The proof is an induction on $k$, the case $k = 0$ being well-known (see for instance \cite[Theorem 19.1.6]{katokhasselblatt}). Let $1 \leq k < r-1$ and assume that $E_u$ is $\mathcal{C}^{\alpha,k-1}$. We will prove that $E_u$ is in $\mathcal{C}^{\alpha_0,k}$ for some $0 < \alpha_0 < 1$.

Let $\mathcal{Q}$ be the section of $\End(TM)$ such that for each $x \in M$ the operator $\mathcal{Q}(x)$ is the orthogonal projection from $T_x M$ to $E_u(x)$, and let $\mathcal{P} = I - \mathcal{Q}$. We want to prove that $\mathcal{Q} \in \mathcal{C}^{\alpha_0,k}\p{M;\End(TM)}$ for some $0 < \alpha_0 \leq \alpha$. Since $E_u$ is $\mathcal{C}^{\alpha,k-1}$ by assumption, we may use any covariant derivative to define the space $\mathcal{C}^{\alpha_0,k}\p{M;\End(TM)}$. We will use the covariant derivative from Remark \ref{remark:different_connection} in order to be able to use \eqref{eq:Leibniz_endomorphism}. Notice that if $X \in \mathcal{C}^{\alpha,k-1}(M;TM)$ is tangent to $E_u$, then $\nabla_X \mathcal{Q}$ is well-defined and continuous by the unstable manifold theorem \cite[\S 3.6]{yoccozIntroductionHyperbolicDynamics1995}, our goal is to prove that it is actually $\mathcal{C}^{\alpha,k-1}$.

We define then a $\mathcal{C}^r$ lift $\Phi : \End(TM) \to \End(TM)$ of $\phi$ by 
\begin{equation*}
\Phi(A) =  \mathrm{d}_x \phi A {}^t \mathrm{d}_x \phi  \textup{ for } A \in \End(T_xM), 
\end{equation*}
where the adjoint is defined using our choice of Riemannian metric on $M$. We write ${}^t B$ for the adjoint of $B$. Let $m \geq 1$ be a large integer and define
\begin{equation*}
V_m(x) = \mathrm{d}_{\phi^{-m}x} \phi^m \mathcal{Q}(\phi^{-m} x) {}^t \mathrm{d}_{\phi^{-m} x} \phi^m  
\end{equation*}
and
\begin{equation*}
W_m(x) = \mathcal{Q}(x) {}^t (\mathrm{d}_{\phi^{-m}x} \phi^m)^{-1} \mathcal{Q}(\phi^{-m}x)(\mathrm{d}_{\phi^{-m} x} \phi^m)^{-1} \mathcal{Q}(x).
\end{equation*}
Notice that $V_m$ and $W_m$ have the same regularity as $\mathcal{Q}$, that is they belong to $\mathcal{C}^{\alpha,k-1}\p{M;\End(TM)}$, and that we have $W_m V_m= V_m W_m = \mathcal{Q}$. Indeed
\begin{equation*}
\begin{split}
& W_m(x) V_m(x)\\ &  = \mathcal{Q}(x) {}^t (\mathrm{d}_{\phi^{-m}x} \phi^m)^{-1} \mathcal{Q}(\phi^{-m}x)(\mathrm{d}_{\phi^{-m} x} \phi^m)^{-1} \mathcal{Q}(x) \\ & \qquad \qquad \qquad \qquad \qquad \qquad \qquad \quad \times \mathrm{d}_{\phi^{-m}x} \phi^m \mathcal{Q}(\phi^{-m} x) {}^t \mathrm{d}_{\phi^{-m} x} \phi^m \\
&  = \mathcal{Q}(x) {}^t (\mathrm{d}_{\phi^{-m}x} \phi^m)^{-1} (\mathrm{d}_{\phi^{-m} x} \phi^m)^{-1} \mathrm{d}_{\phi^{-m}x} \phi^m \mathcal{Q}(\phi^{-m} x) {}^t \mathrm{d}_{\phi^{-m} x} \phi^m \\
& = \mathcal{Q}(x) {}^t (\mathrm{d}_{\phi^{-m}x} \phi^m)^{-1} \mathcal{Q}(\phi^{-m} x) {}^t \mathrm{d}_{\phi^{-m} x} \phi^m \\
& = {}^t\p{\mathrm{d}_{\phi^{-m} x} \phi^m \mathcal{Q}(\phi^{-m} x) (\mathrm{d}_{\phi^{-m}x} \phi^m)^{-1} \mathcal{Q}(x)}\\ & = {}^t\p{\mathrm{d}_{\phi^{-m} x} \phi^m (\mathrm{d}_{\phi^{-m}x} \phi^m)^{-1} \mathcal{Q}(x)} \\
& = {}^t \mathcal{Q}(x) = \mathcal{Q}(x).
\end{split}
\end{equation*}
Here we used on the third and sixth line that the derivative of $\phi$ preserves the unstable direction. We get the other equality from $V_m W_m = {}^t\p{W_m V_m} = {}^t \mathcal{Q} = \mathcal{Q}$. We also have $\mathcal{L}_{\Phi}^m(\mathcal{Q}) = V_m$. Applying Lemma \ref{lemma:chain_rule}, we find, for $X \in \mathcal{C}^{\alpha,k-1}\p{M;TM}$ parallel to the unstable direction,
\begin{equation*}
\begin{split}
\nabla_X \mathcal{Q} & = \nabla_X(W_m V_m) \\ & = (\nabla_X W_m) V_m + W_m \mathcal{L}_{\Phi}^m(\nabla_{\phi^{-m}_* X} \mathcal{Q})+ W_m \gamma_m(\mathcal{Q} \circ \phi^{-m}, X),
\end{split}
\end{equation*}
for some $\mathcal{C}^{\alpha,k-1}$ section $\gamma_m$ of $(\phi^{-m})^*\End(TM)^* \otimes T^* M \otimes \End(TM)$ (actually $\gamma_m$ may be shown to be $\mathcal{C}^{r-1}$). Since ${}^t \mathrm{d} \phi^m$ preserves the orthogonal of the unstable direction, we have that $V_m \mathcal{P}= 0$, which implies
\begin{equation*}
(\nabla_X \mathcal{Q}) \mathcal{P} = W_m \mathcal{L}_{\Phi}^m(\nabla_{\phi^{-m}_* X} \mathcal{Q}) \mathcal{P} + W_m \gamma_m(\mathcal{Q} \circ \phi^{-m}, X) \mathcal{P}.
\end{equation*}
Similarly, we find that
\begin{equation*}
\begin{split}
\nabla_X \mathcal{Q} & = \nabla_X(V_m W_m) \\ & = V_m (\nabla_X W_m) + \mathcal{L}_{\Phi}^m(\nabla_{\phi^{-m}_* X} \mathcal{Q})W_m+  \gamma_m(\mathcal{Q} \circ \phi^{-m}, X) W_m,
\end{split}
\end{equation*}
so that
\begin{equation*}
\mathcal{P} (\nabla_X \mathcal{Q}) = \mathcal{P} \mathcal{L}_{\Phi}^m(\nabla_{\phi^{-m}_* X} \mathcal{Q})W_m+  \mathcal{P} \gamma_m(\mathcal{Q} \circ \phi^{-m}, X) W_m.
\end{equation*}
Differentiating $\mathcal{Q}^2  = \mathcal{Q}$, we find $(\nabla_X \mathcal{Q}) \mathcal{Q} = \mathcal{P} (\nabla_X \mathcal{Q})$, so that $\nabla_X \mathcal{Q} = \mathcal{P} \nabla_X \mathcal{Q} + (\nabla_X \mathcal{Q}) \mathcal{P}$, and thus,
\begin{equation}\label{eq:la_relation}
\begin{split}
\nabla_X \mathcal{Q} & = W_m \mathcal{L}_{\Phi}^m(\nabla_{\phi^{-m}_* X} \mathcal{Q}) \mathcal{P} + \mathcal{P} \mathcal{L}_{\Phi}^m(\nabla_{\phi^{-m}_* X} \mathcal{Q})W_m \\ & \qquad \qquad \qquad + W_m \gamma_m(\mathcal{Q} \circ \phi^{-m}, X) \mathcal{P} + \mathcal{P} \gamma_m(\mathcal{Q} \circ \phi^{-m}, X) W_m.
\end{split}
\end{equation}
Let then $X_1,\dots,X_D$ be as in the definition of the norm \eqref{eq:def_norm} and the $a_{ij}^m$'s be as in \eqref{eq:decompostion_Xi}. Consider the section $\mathcal{R} = (\mathcal{R}_1,\dots,\mathcal{R}_D)$ of $\End(TM)^{\oplus D}$ defined by $\mathcal{R}_i = \nabla_{X_i} \mathcal{Q}$. Using \eqref{eq:la_relation} and \eqref{eq:decompostion_Xi}, we find the that for $i = 1,\dots,D$, we have
\begin{equation}\label{eq:equation_derivative}
\begin{split}
\mathcal{R}_i & = \sum_{j =1}^D a_{ij}^m \circ \phi^{-m} (W_m \mathcal{L}_{\Phi}^m(\mathcal{R}_j) \mathcal{P} + \mathcal{P}\mathcal{L}_{\Phi}^m(\mathcal{R}_j) W_m) \\ & \qquad \qquad \qquad + W_m \gamma_m(\mathcal{Q} \circ \phi^{-m}, X_i) \mathcal{P} + \mathcal{P} \gamma_m(\mathcal{Q} \circ \phi^{-m}, X_i) W_m.
\end{split}
\end{equation}
As in the proof of the general case in Proposition \ref{prop:bound_spectral_radius}, it follows from \eqref{eq:equation_derivative} that there is a lift $\hat{\Phi}_m : \End(TM)^{\oplus D} \to \End(TM)^{\oplus D}$ such that
\begin{equation}\label{eq:characterization_R}
\mathcal{R} = \mathcal{L}_{\hat{\Phi}_m}(\mathcal{R}) + \mathcal{H},
\end{equation}
where $\mathcal{H} = (\mathcal{H}_1,\dots,\mathcal{H}_D)$ with $\mathcal{H}_i = W_m \gamma_m(\mathcal{Q} \circ \phi^{-m}, X_i) \mathcal{P} + \mathcal{P} \gamma_m(\mathcal{Q} \circ \phi^{-m}, X_i) W_m$ for $i = 1,\dots,D$. It follows from our induction hypothesis that $\hat{\Phi}_m$ and $\mathcal{H}$ are $\mathcal{C}^{\alpha,k-1}$. Let us prove that $\theta(\hat{\Phi}_m) < 1$ when $m$ is large enough. We can use \eqref{eq:borne_aij} to bound the $a_{ij}^m$'s. Hence, we only need to prove that if $x \in M$ and $A \in \End(T_x M)$ has norm $1$, then the norm of 
\begin{equation}\label{eq:bout_de_lift}
W_m(\phi^m x) \Phi^m(A) \mathcal{P}(\phi^m x) + \mathcal{P}(\phi^m x) \Phi^m(A) W_m(\phi_m x)
\end{equation}
is bounded uniformly in $m$. To do so, we start by writing
\begin{equation}\label{eq:un_calcul}
\begin{split}
& W_m(\phi^m x) \Phi^m(A) \mathcal{P}(\phi^m x) \\ & \qquad = \mathcal{Q}(\phi^m x) {}^t (\mathrm{d}_{x} \phi^m)^{-1} \mathcal{Q}(x)(\mathrm{d}_{ x} \phi^m)^{-1} \mathcal{Q}( \phi^m x) \mathrm{d}_{x} \phi^m A {}^t \mathrm{d}_{x} \phi^m \mathcal{P}(\phi^m x).
\end{split}
\end{equation}
We see that the operator on the right of $A$ is ${}^t\p{\mathcal{P}(\phi^m x)\mathrm{d}_{x} \phi^m}$, which is exponentially small in $m$ due to the hyperbolicity of $\phi$. Let us then consider the operator on the left of $A$. To do so, introduce the projector $\pi(x)$ on $E_s(x)$ along $E_u(x)$. Using that the unstable direction is stable under the action of the derivative of $\phi$, we find that the operator on the left of $A$ in  \eqref{eq:un_calcul} is
\begin{equation*}
{}^t\p{ (\mathrm{d}_x \phi^m)^{-1} \mathcal{Q}(\phi^m x)} \p{\mathcal{Q}(x)(I - \pi(x)) + \mathcal{Q}(x) (\mathrm{d}_x \phi^m)^{-1} \mathcal{Q}(\phi^m x) \mathrm{d}_x \phi^m \pi(x)},
\end{equation*}
and it follows from the hyperbolicity of $\phi$ that this operator is exponentially small in $m$. Hence, we see that the term on the left in \eqref{eq:bout_de_lift} is exponentially small in $m$. By taking the adjoint, we see that the same is true for the term on the right. Consequently, we have $\theta(\hat{\Phi}_m) < 1$, provided $m$ is large enough.

Choosing $\alpha_0$ small enough so that $\theta(\hat{\Phi}_m) \Lambda^{\alpha_0 m} <1$, we find that the equation \eqref{eq:characterization_R} has a unique solution $(I - \mathcal{L}_{\hat{\Phi}_m})^{-1} \mathcal{H}$ both in $\mathcal{C}^{0,0}(M;\End(TM)^{\oplus D})$ and in $ \mathcal{C}^{\alpha,k-1}(M;\End(TM)^{\oplus D})$. Hence, these solutions coincide and since we know that $\mathcal{R}$ is continuous, it follows that $\mathcal{R}$ is infact in $\mathcal{C}^{\alpha,k-1}$. Thus, for every $X \in \mathcal{C}^{\alpha,k-1}(M;TM)$ parallel  to $E_u$, we have that $\nabla_X \mathcal{Q}$ is $\mathcal{C}^{\alpha,k-1}$ and it follows that $\mathcal{Q}$ is $\mathcal{C}^{\alpha,k}$, ending the proof of the lemma.
\end{proof}

\begin{remark}\label{remark:unstable_basis}
Notice that a consequence of the proof of Lemma \ref{lemma:regularity} is that one can compute $\nabla_X \mathcal{Q}$ where $X$ is a vector field parallel to the unstable direction (i.e. one can compute ``the derivatives of $E_u$ in the unstable direction''). Indeed, using the notation from the proof of Lemma \ref{lemma:regularity}, one can use the method from Remark \ref{remark:solving_equations} to deduce $\mathcal{R}$ from \eqref{eq:characterization_R}. 

As a consequence, one can construct a family $X_1,\dots,X_D$ of $\mathcal{C}^{k,\alpha}$ vector fields (where $k$ and $\alpha$ are as in Lemma \ref{lemma:regularity}), parallel to the unstable direction such that, for every $x \in M$, the vectors $X_1(x),\dots,X_D(x)$ span $E_u(x)$ and for which one can compute $\nabla_X X_1,\dots,\nabla_X X_D$ when $X$ is parallel to the unstable direction. Let us explain how. One just needs to do it locally (thanks to the compactness of $M$). Let $x_0 \in M$ and choose $\widetilde{X}_1,\dots,\widetilde{X}_{d_u}$ a family of smooth vector fields such that $\widetilde{X}_1(x_0),\dots,\widetilde{X}_{d_u}(x_0)$ is a basis for $E_u(x_0)$. Then, set $X_j(x) = \mathcal{Q}(x) \widetilde{X}_j(x)$ for $j = 1,\dots,d_u$ and $x \in M$. Notice that for $x$ near $x_0$ the family $X_1(x),\dots,X_{d_u}(x)$ is a basis for $E_u(x)$. Finally, if $X$ is parallel to the unstable direction, we have
\begin{equation*}
\nabla_X X_j = (\nabla_X Q) \widetilde{X}_j + Q (\nabla_X \widetilde{X}_j).
\end{equation*}
\end{remark}

\section{S3 algorithm decomposition}\label{section:Y_and_Z}

In this section, we explain how the vector fields $Y$ and $V$ are constructed in the second step of the S3 algorithm. Here, we apply Proposition \ref{prop:bound_spectral_radius} to get a statement that encompasses \cite[Lemma 8.6 and Proposition 8.9]{chandramoorthyEfficientComputationLinear2021}.

\begin{prop}\label{prop:decomposition}
There is $\alpha_0 > 0$ such that, if $0 \leq \alpha \leq \alpha_0,k < r-1$ and $X \in \mathcal{C}^{\alpha,k}\p{M;TM}$, then there are $Y$ and $V$ in $\mathcal{C}^{\alpha,k}\p{M;TM}$ such that $Y$ is parallel to the unstable direction and
\begin{equation*}
X = Y + V - \phi_*V.
\end{equation*}
\end{prop}

In order to apply Proposition \ref{prop:bound_spectral_radius} to prove Proposition \ref{prop:decomposition}, we need the following estimate.

\begin{lemma}\label{lemma:bound_PdT}
Let $\mathcal{P} : TM \to TM$ be the map that acts on each fiber by orthogonal projection on the orthogonal of $E_u$. Then $\mathcal{P} \mathrm{d}\phi : TM \to TM$ is a lift of $\phi$ and $\theta(\mathcal{P} \mathrm{d}\phi) <1$.
\end{lemma}

\begin{proof}
Let $\pi : TM \to TM$ acts on each fiber by projection on $E_s$ along $E_u$. Since $\pi$ and $\mathcal{P}$ are projectors with the same kernel in each fiber, we have that
\begin{equation}\label{eq:bound_PdT}
\mathcal{P} \pi = \mathcal{P} \quad \textup{and} \quad \pi \mathcal{P} = \pi.
\end{equation}
We also have $\mathcal{\pi} \mathrm{d}\phi = \mathrm{d}\phi \pi$. Hence, for $n \geq 1$, we have
\begin{equation*}
\begin{split}
\p{\mathcal{P} \mathrm{d}\phi}^n & = \p{\mathcal{P} \pi \mathrm{d}\phi}^n = \p{\mathcal{P} \mathrm{d}\phi \pi}^n \\
       & = \mathcal{P} \p{\mathrm{d}\phi \pi \mathcal{P}}^{n-1} \mathrm{d}\phi \pi = \mathcal{P} \p{\mathrm{d}\phi \pi}^{n}.
\end{split}
\end{equation*}
From the hyperbolicity of $\phi$, we know that operator norm of $\p{\mathrm{d}\phi \pi}^n$ in the fibers decays exponentially fast with $n$, so that we have $\theta\p{\mathcal{P}\mathrm{d}\phi} < 1$.
\end{proof}

\begin{proof}[Proof of Proposition \ref{prop:decomposition}]
Let $\Phi = \mathcal{P} \mathrm{d}\phi$ be the lift of $\phi$ from Lemma \ref{lemma:bound_PdT}. It follows from Lemma \ref{lemma:regularity} that $\Phi$ is $\mathcal{C}^{\alpha,k}$ for $\alpha > 0$ small. According to Lemma \ref{lemma:bound_PdT}, we have $\theta(\Phi) < 1$.

Let then $\alpha_0$ be such that $\theta(\Phi) \Lambda^{\alpha_0} < 1$ and $E_u$ is $\mathcal{C}^{\alpha_0,k}$. Then, if $\alpha \leq \alpha_0$, the spectral radius of $\mathcal{L}_{\Phi}$ acting on $\mathcal{C}^{\alpha,k}\p{M;TM}$ is strictly less than $1$. Consequently, if $X \in \mathcal{C}^{\alpha,k}\p{M;TM}$, we can set
\begin{equation}\label{eq:def_V}
V = \p{I - \mathcal{L}_{\Phi}}^{-1} X \quad \textup{and} \quad Y = X - \p{V - \phi_*V}.
\end{equation}
It follows then from the definition of $\Phi$ that $Y$ is parallel to the unstable direction. Indeed, recalling the orthogonal projector $\mathcal{Q}$ on $E_u$, we have
\begin{equation*}
Y = X - (V - \phi_*V) = X - (V - \mathcal{P}\phi_* V) + \mathcal{Q}\phi_* V = \mathcal{Q}\phi_* V.
\end{equation*}
\end{proof}

\begin{remark}\label{remark:step_2}
Proposition \ref{prop:decomposition} implies that there exist vector fields $Y$ and $V$ as in the second step of the S3 algorithm. Moreover, the proof of Proposition \ref{prop:decomposition} gives a way to compute $Y$ and $V$. Indeed, if we look at the expression \eqref{eq:def_V} for $V$ and $Y$, we see that their values at a point $x \in M$ may be approximated from the knowledge of $X$ and $E_u$ on the backward orbit of $x$, as explained in Remark \ref{remark:solving_equations}. In addition, the quality of the approximation will be exponentially small in the number of points in the backward orbit of $x$ that we use. Here, we recall that $X$ is given and that $E_u$ has been computed in the first step of the algorithm.
\end{remark}

\begin{remark}
Notice that in Proposition \ref{prop:decomposition}, even if $X$ is $\mathcal{C}^{k+ \alpha}$, then $Y$ and $V$ will \emph{a priori} only be $\mathcal{C}^{\alpha,k}$. This is due to the phenomenon highlighted in Example \ref{example:basic} that causes loss of regularity in the stable direction. However, this is not a problem since, even if we cannot work with mere Hölder regularity (because of the use of integration by parts for instance), having a Hölder-continuous derivative in the unstable direction is enough for the arguments from \S \ref{section:unstable} and \S \ref{section:convergence}.
\end{remark}

\begin{remark}
It is not crucial in the proof of Proposition \ref{prop:decomposition} that $\mathcal{P}$ is an orthogonal projector. The only things that matter is that $\mathcal{P}$ is a projector with kernel $E_u$ (so that \eqref{eq:bound_PdT} holds) and that $\mathcal{P}$ is $\mathcal{C}^{\alpha,k}$ in order to get the same regularity for $V$ and $Y$. Notice for instance that another choice of Riemannian metric would \emph{a priori} leads to another projector $\mathcal{P}$ and hence to another choice of $Y$ and $V$ in Proposition \ref{prop:decomposition}.
\end{remark}

\section{Unstable contribution}\label{section:unstable}

In this section, we give a method to compute the quantity $\Psi_{\phi}(Y,f)$, defined by \eqref{eq:result}, when $Y$ is parallel to the unstable direction, justifying the third step in the S3 algorithm. The language from differential geometry will allow us to deal with this unstable contribution without relying on a particular parametrization of the unstable manifold (as it is the case in \cite{chandramoorthyEfficientComputationLinear2021,sliwiak2021spacesplit}).

\begin{prop}\label{prop:adjoint_Y}
Let $0 \leq \alpha  <1$ and $1 \leq k < r - 1$ be such that $E_u$ is $\mathcal{C}^{\alpha,k}$. Let $Y \in \mathcal{C}^{\alpha,k}\p{M;TM}$ be parallel to the unstable direction. Then, there is a continuous function $\rho_Y : M \to \C$ such that if we set $Y^* = - Y + \rho_Y$ then, for every  $\mathcal{C}^{1+\alpha}$ functions $f,g$, we have
\begin{equation*}
\mu\p{Yf. g} = \mu\p{f.Y^* g}.
\end{equation*}
\end{prop}

\begin{proof}
Cover $M$ by open sets $\p{U_\omega}_{\omega \in \Omega}$ foliated by local unstable manifolds. Then, let $\p{h_\omega}_{\omega \in \Omega}$ be a partition of unity subordinated to the open cover $\p{U_\omega}_{\omega \in \Omega}$. We can then split
\begin{equation*}
\mu\p{Yf. g} = \sum_{\omega \in \Omega}\mu\p{Yf. h_\omega g}
\end{equation*}
and consider each term separately. Fix then $\omega \in \Omega$ and let $W_\omega$ denote a transverse to the unstable foliation in $U_\omega$. We can then disintegrate the restriction of $\mu$ to $U_\omega$ and write
\begin{equation*}
\begin{split}
\mu\p{Yf. h_\omega g} = \int_{W_\omega} \p{\int_{W_{u,\textup{loc}}(y)} Yf(x) h_\omega (x) g(x) \mathrm{d}m_y(x) } \mathrm{d}\nu_\omega(y),
\end{split}
\end{equation*}
where $\nu_\omega$ is a measure on $W_\omega$ and, for every $y \in W_\omega$, the measure $m_y$ is absolutely continuous and fully supported on the local unstable manifold $W_{u,\textup{loc}}(y)$. Consequently, for $y \in W_\omega$, we can integrate by parts in the integral on $W_{u,\textup{loc}}(y)$ to find
\begin{equation*}
\begin{split}
& \int_{W_{u,\textup{loc}}(y)} Yf(x) h_\omega (x) g(x) \mathrm{d}m_y(x) \\ & \qquad \qquad = \int_{W_{u,\textup{loc}}(y)} f(x) \p{ - h_\omega (x) Y g(x) + \rho_\omega(x) g(x)} \mathrm{d}m_y(x),
\end{split}
\end{equation*}
for some bounded measurable function $\rho_\omega$. Here, the integration by part is licit since the measure $\mathrm{d}m_y$ has a differentiable density and the restriction of $Y$ to an unstable manifold is differentiable, and there is no boundary term since $h_\omega$ vanishes near the boundary of the local unstable manifold. Moreover, using the explicit formula for the density of $\mathrm{d}m_y$ given for instance in \cite[Corollary 6.1.4]{ledrappierMetricEntropyDiffeomorphisms1985a} and the fact that the derivatives of $Y$ in the unstable direction are Hölder, we see that $\rho_\omega$ is actually continuous. Summing over $\omega$, the proposition is proven. 
\end{proof}

\begin{lemma}\label{lemma:formula_rho}
Let $0 \leq \alpha  <1$ and $1 \leq k < r - 1$ be such that $E_u$ is $\mathcal{C}^{\alpha,k}$. Let $X_1,\dots,X_D$ be as in \eqref{eq:def_norm}. Let the $a_{ij}^m$'s be as in \eqref{eq:decompostion_Xi}. Then, for every $m \geq 1$ and $i = 1,\dots,D$, we have
\begin{equation}\label{eq:formula_rho}
\rho_{X_i} = \sum_{j = 1}^D a_{ij}^m \circ  \phi^{-m} \rho_{X_j} \circ \phi^{-m} - \sum_{j = 1}^D X_j a_{ij}^m \circ \phi^{-m}.
\end{equation}
Let $Y \in \mathcal{C}^{\alpha,k}\p{M;TM}$ be parallel to the unstable direction and use Lemma \ref{lemma:spaces} to write $Y = \sum_{j = 1}^D b_j X_j$, where the coefficients $b_j$ are $\mathcal{C}^{\alpha,k}$, then
\begin{equation}\label{eq:general_rho_Y}
\rho_Y = \sum_{j = 1}^D b_j \rho_{X_j} - X_j b_j.
\end{equation}

In particular, there is $\beta >0$ such that if $Y \in \mathcal{C}^{\alpha,k}\p{M;TM}$ is parallel to the unstable direction, then $\rho_Y$ is $\mathcal{C}^{\beta,k-1}$.
\end{lemma}

\begin{proof}
For $f$ smooth, we have
\begin{equation*}
\begin{split}
& \mu\p{(\sum_{j = 1}^D a_{ij}^m \circ \phi^{-m} \rho_{X_j} \circ \phi^{-m} - \sum_{j = 1}^D X_j a_{ij}^m \circ \phi^{-m})f} \\ & \qquad \qquad \qquad = \mu\p{(\sum_{j = 1}^D a_{ij}^m \rho_{X_j} - \sum_{j = 1}^D X_j a_{ij}^m)f \circ \phi^m} \\
   & \qquad \qquad \qquad = \mu\p{\sum_{j =1}^D a_{ij}^m X_j(f \circ \phi^m)} = \mu\p{ (\phi^{-m}_* X_i)(f \circ \phi^m)}\\ 
   & \qquad \qquad \qquad = \mu(X_i f \circ \phi^m) =\mu(X_i f) = \mu(\rho_{X_i} f),
\end{split}
\end{equation*}
and \eqref{eq:formula_rho} follows since $\mu$ has full support. Let us now show that $\rho_{X_1},\dots,\rho_{X_D}$ are $\mathcal{C}^{\beta,k-1}$ for some $\beta > 0$. To do so let $\boldsymbol{\rho} = (\rho_{X_1},\dots,\rho_{X_D})$. Let $m \geq 1$ be a large integer and let $w^m = (w_1^m,\dots,w_D^m)$ be the section of $M \times \mathbb{C}^D$ defined by
\begin{equation*}
w_i^m = - \sum_{j = 1}^D X_j a_{ij}^m \circ \phi^{-m} \textup{ for } i = 1,\dots,D.
\end{equation*}
Let also $\Phi_m$ be the lift of $\phi^m$ on $M \times \mathbb{C}^D$ given by the formula
\begin{equation*}
\Phi_m( x, (z_i)_{1 \leq i \leq D}) = (\phi^m(x), (\sum_{j = 1}^D a_{ij}^m z_j)_{1 \leq i \leq D}).
\end{equation*}
It follows then from \eqref{eq:formula_rho} that we have
\begin{equation}\label{eq:characterization_rho}
\p{I - \mathcal{L}_{\Phi_m}} \boldsymbol{\rho} = w.
\end{equation}
Moreover, it follows from \eqref{eq:borne_aij} that, for $m$ large enough, we have $\theta(\Phi_m) < 1$. Hence, since $\boldsymbol{\rho}$ is continuous, it follows from \eqref{eq:characterization_rho} and Proposition \ref{prop:bound_spectral_radius} that $\boldsymbol{\rho}$ is $\mathcal{C}^{\beta,k-1}$ if $\beta$ is small enough so that $\theta(\Phi_m)\Lambda^{\beta m} < 1$ and $w$ is $\mathcal{C}^{\beta,k-1}$ itself.

Finally, to deal with the case of a general $Y \in \mathcal{C}^{\alpha,k}\p{M;TM}$ parallel to the unstable direction, just notice that for $f$ smooth we have
\begin{equation*}
\begin{split}
\mu(f \rho_Y)& = \mu(Yf) = \sum_{j = 1}^D \mu(X_jf b_j) = \sum_{j = 1}^D \mu(f(-X_j b_j + b_j \rho_{X_j})) \\ & = \mu\p{f \p{\sum_{j = 1}^D b_j \rho_{X_j} - X_j b_j}}.
\end{split}
\end{equation*}
\end{proof}

\begin{remark}\label{remark:compute_rho}
As explained in the introduction, the knowledge of $\rho_Y$ allows us to approximate the contribution $\Psi_{\phi}(Y,f)$ of $Y$ to the linear response using the formula \eqref{eq:unstable_contribution}. Notice that the characterization of $\rho_Y$ given in Lemma \ref{lemma:formula_rho} enables to use the method explained in Remark \ref{remark:solving_equations} to compute $\rho_Y$ from the knowledge of $E_u,Y$ and $\nabla_{X} Y$ for $X$ parallel to the unstable direction. Indeed, it follows from \eqref{eq:general_rho_Y} that we need to compute the coefficients $b_j$, their derivatives $X_j b_j$ and $\rho_{X_1},\dots,\rho_{X_D}$. To compute $\rho_{X_1},\dots,\rho_{X_D}$ using \eqref{eq:formula_rho} (or its rephrasing \eqref{eq:characterization_rho}) and Remark \ref{remark:solving_equations}, we need to compute the coefficients $a_{ij}^m$ and their derivatives $X_j a_{ij}^m$. To compute the $b_j, X_j b_j, a_{ij}^m$ and $X_j a_{ij}^m$, we can use the strategy from Remark \ref{remark:coefficients}, which requires to know $Y,X_1,\dots,X_D$ and the derivatives $\nabla_{X_j} Y$ and $\nabla_{X_j} X_i$ for $1 \leq i,j \leq D$. The derivative $\nabla_{X_j} Y$ is computed using the strategy from Remark \ref{remark:solving_derivative} and the definition \eqref{eq:def_V} of $Y$. To compute the $\nabla_{X_j} X_i$, one can rely on Remark \ref{remark:unstable_basis}.

As usual, the value of $\rho_Y$ at $x$ is obtained with a precision which is exponential in the number of points in the backward orbit of $x$ that we use for the computation.
\end{remark}

\section{Coboundary contribution}\label{section:convergence}

In order to end the proof of Theorem \ref{theorem:main}, we want now to evaluate the contribution to $\Psi_{\phi}(X,f)$ of the coboundary term in the decomposition \eqref{eq:decomposition_X}. As mentioned in the introduction, we depart from the proof strategy in \cite{chandramoorthyEfficientComputationLinear2021}. In particular, here we use the spaces of anisotropic distributions from \cite{GLK}, that allows us to work directly with the term $\Psi_\phi(V - \phi_*V,f)$ in \eqref{eq:decomposition_X} (the approach in \cite{chandramoorthyEfficientComputationLinear2021} uses the intermediates in the computation of $V$). The main result from this section is:

\begin{lemma}\label{lemma:telescopic}
Let $0 < \alpha <1$ be such that $E_u$ is $\mathcal{C}^\alpha$, a vector field $V \in \mathcal{C}^{\alpha,1}\p{M;TM}$ and $f$ be a $\mathcal{C}^{1+\alpha}$ function on $M$. 
Then the series defining $\Psi_{\phi}\p{V - \phi_*V,f}$ converges and 
\begin{equation*}
\Psi_{\phi}\p{V - \phi_*V,f} = \mu\p{Vf}.
\end{equation*}
\end{lemma}

Lemma \ref{lemma:telescopic} follows from the estimate:

\begin{lemma}\label{lemma:convergence}
Let $0 <\alpha<1$ be such that $E_u$ is $\mathcal{C}^\alpha$, a vector field $V \in  \mathcal{C}^{\alpha,1}\p{M;TM}$ and $f \in \mathcal{C}^{1+\alpha}$. Then the integral $\mu\p{V(f \circ \phi^n)}$ tends to $0$ (exponentially fast) as $n$ tends to $+ \infty$.
\end{lemma}

\begin{proof}[Proof of Lemma \ref{lemma:telescopic}]
Just notice that for $n \in \N$, we have
\begin{equation*}
\begin{split}
\mu\p{\p{V - \phi_* V}\p{f \circ \phi^n}} & = \mu\p{V \p{ f \circ \phi^n}} - \mu\p{V\p{f \circ \phi^{n+1}} \circ \phi^{-1}} \\
     & = \mu\p{V \p{ f \circ \phi^n}} - \mu\p{V\p{f \circ \phi^{n+1}}}.
\end{split}
\end{equation*}
Hence, $\Psi_{\phi}(V - \phi_* V,f)$ becomes a telescopic sum, and we get thanks to Lemma \ref{lemma:convergence}:
\begin{equation*}
\Psi_{\phi}(V - \phi_* V,f) = \mu\p{V f}.
\end{equation*}
\end{proof}

The rest of this section is dedicated to the proof of Lemma \ref{lemma:convergence}. We will use the machinery from \cite{GLK} (see also \cite{GouLiv2}). The spaces of anisotropic distribution introduced in this reference allows us to deal directly with the coboundary contribution $\Psi_\phi(V - \phi_* V,f)$, while the approach from \cite{chandramoorthyEfficientComputationLinear2021} is more indirect.

Let us remind a few facts from \cite{GLK}. Letting $0 <  q <\min(\alpha,r-2)$, we recall the space $\mathcal{B}^{1,q}$ from \cite{GLK}. For convenience, we will reverse time with respect to \cite{GLK,GouLiv2} and consider the Koopman operator
\begin{equation}\label{eq:koopman_operator}
\mathcal{K} : u \mapsto u \circ \phi
\end{equation}
instead of the transfer operator. To work with the Koopman operator, we only need to replace the admissible stable leaves from \cite{GLK} with admissible unstable leaves (which, in \cite{GLK}, would be admissible leaves for $\phi^{-1}$)  in the definition of the space $\mathcal{B}^{1,q}$. It amounts to replace $\phi$ by $\phi^{-1}$ and consider a weighted transfer operator. This case is not exactly included in the analysis from \cite{GLK}, but very few changes are needed to adapt the results from \cite{GLK} to the operator \eqref{eq:koopman_operator}, as the interested reader can check it (the main change would be in \cite[Lemma 6.2]{GLK}, whose proof actually becomes simpler since there is no Jacobian involved). The technical tools to deal with much more general weighted transfer operators are exposed in \cite{GouLiv2}.

The space $\mathcal{B}^{1,q}$ is defined as the completion of $\mathcal{C}^{1 + \alpha}$ with respect to the norm 
\begin{equation*}
\n{u}_{\gl} \coloneqq \max\p{\va{u}_{0,q},\va{u}_{1,1+q}},
\end{equation*}
where
\begin{equation*}
\va{u}_{0,q} = \sup_{W \in \Sigma} \sup_{\substack{g \in \mathcal{C}_0^q\p{W,\R} \\ \va{g}_{\mathcal{C}^q} \leq 1}} \int_W u g \mathrm{d}m_{W}
\end{equation*}
and
\begin{equation*}
\va{u}_{1,1+q} = \sup_{W \in \Sigma} \sup_{\substack{v \in \mathcal{V}^{1 + \alpha}(W) \\ \va{v}_{\mathcal{C}^{1 + \alpha}} \leq 1}} \sup_{\substack{g \in \mathcal{C}_0^{1+q}\p{W,\R} \\ \va{g}_{\mathcal{C}^{1+q}} \leq 1}} \int_W (vu) g \mathrm{d}m_W.
\end{equation*}
Here, $\Sigma$ denotes a set of admissible unstable leaves (that is a set of admissible leaves, as defined in \cite{GLK}, but for $\phi^{-1}$). We will only need to know that the elements of $\Sigma$ are small $d_u$-dimensional disks in $M$ and that there is $\delta > 0$ such that any disk of radius $\delta$ in an unstable manifold is contained in an element of $\Sigma$, which is itself a piece of unstable manifold (see \cite[Definition 3.1]{GouLiv2}). For $W \in \Sigma$, we write $\mathrm{d}m_W$ for the Lebesgue measure on $W$, we denote by $\mathcal{V}^{1 + \alpha}\p{W}$ the set of $\mathcal{C}^{1 + \alpha}$ vector fields on a neighbhourhood of $W$ and by $\mathcal{C}_0^q\p{W}$ the space of $\mathcal{C}^q$ functions on $W$ supported away from the boundary. Notice that we recall this definition of the space $\mathcal{B}^{1,q}$ mostly for educational purpose, as we will only need the local description of the space $\mathcal{B}^{1,q}$ given by \cite[Lemma 3.2]{GLK}.

From \cite[Theorem 2.3]{GLK}, or rather its adaption to the case of Koopman operator instead of transfer operator, we know that the essential spectral radius of $\mathcal{K}$ acting on $\mathcal{B}^{1,q}$ is strictly less than $1$, and that $1$ is a simple eigenvalue for $\mathcal{K}$. Moreover, this is the only eigenvalue of $\mathcal{K}$ on the unit circle, the associated right eigenspace consists of constant functions, while the left eigenspace is generated by the SRB measure $\mu$ of $T$. Consequently, we may decompose $\mathcal{K}$ as
\begin{equation}\label{eq:spectral_decomposition}
\mathcal{K} = P +Q,
\end{equation}
where $P$ is the rank $1$ operator defined by $P(u) = \mu\p{u}$ (here we identify a number with a constant function) and $Q$ is a bounded operator on $\mathcal{B}^{1,q}$ with spectral radius strictly less than $1$ that satisfies $PQ = QP = 0$. Consequently, for $n \geq 1$, we have
\begin{equation}\label{eq:iterated_spectral_decomposition}
\mathcal{K}^{n} = P + Q^n.
\end{equation}

With these tools from \cite{GLK}, Lemma \ref{lemma:convergence} follows from:

\begin{lemma}\label{lemma:linear_form}
Let $0 < \alpha < 1$ be such that $E_u$ is $\mathcal{C}^\alpha$ and $V \in \mathcal{C}^{\alpha,1}\p{M;TM}$. Then the map
\begin{equation*}
f \mapsto \mu\p{Vf}
\end{equation*}
extends to a bounded linear form on $\mathcal{B}^{1,q}$.
\end{lemma}

Let us explain how Lemma \ref{lemma:linear_form} allows us to prove Lemma \ref{lemma:convergence}.

\begin{proof}[Proof of Lemma \ref{lemma:convergence}]
Since $f \in \mathcal{C}^{1+ \alpha}$, we see that $f$ belongs to $\mathcal{B}^{1,q}$. Consequently, we have for $n \in \N$
\begin{equation*}
f \circ \phi^n = \mathcal{K}^n f = \mu(f) + Q^n f
\end{equation*}
and thus
\begin{equation*}
\mu\p{V(f \circ \phi^n)} = \mu\p{V(f \circ \phi^n - \mu(f))} = \mu\p{V Q^n f}.
\end{equation*}
Here, the expression $\mu\p{V Q^n f}$ makes sense due to Lemma \ref{lemma:linear_form}, and we have for some $C >0$ that does not depend on $n$:
\begin{equation*}
\va{\mu\p{V Q^n f}} \leq C \n{Q^n f}_{\gl}.
\end{equation*}
However, since the spectral radius of $Q$ is strictly less than $1$, we see that $\n{Q^n f}_{\gl}$ decays exponentially fast when $n$ tends to $+ \infty$.
\end{proof}

We need now to prove Lemma \ref{lemma:linear_form}.

\begin{proof}[Proof of Lemma \ref{lemma:linear_form}]
Let $f \in \mathcal{C}^{1 + \alpha}$. Let $\chi$ be a $\mathcal{C}^\infty$ function on $M$ with small support. We will estimate $\mu\p{ \chi Vf}$ using the $\mathcal{B}^{1,q}$ norm of $f$ under the assumption that the support of $\chi$ is small enough. The result then follows by a partition of unity argument.

By assuming that the support of $\chi$ is small enough, we may work in coordinates. In these coordinates, we can write
\begin{equation*}
V = \sum_{j = 1}^d v_j \frac{\partial \phantom{x}}{\partial x_j}.
\end{equation*}
Here, the $v_j$'s are $\mathcal{C}^{\alpha,1}$. In particular, the restriction of the $v_j$'s to unstable leaves are $\mathcal{C}^{1 + \alpha}$ with uniform $\mathcal{C}^{1 + \alpha}$ estimates. Then, up to taking the support of $\chi$ smaller, we can cover the patch of coordinates we are working in by a family $\p{W_z}_{z \in Z}$ of local (exact) unstable leaves of uniform size. We endow $Z$ with a measurable structure by identifying it with a manifold transverse to the unstable direction. Moreover, we may assume that the $W_z$'s belong to $\Sigma$ and that each $W_z$ fully intersects the support of $\chi$ (i.e. the support of the restriction of $\chi$ to $W_z$ does not intersect the boundary of $W_z)$. Then, we disintegrate $\mu$ (restricted to the coordinate patch) with respect to the $W_z$'s and find that, for $j =1,\dots,d$, we have
\begin{equation*}
\mu\p{\chi v_j \frac{\partial f}{\partial x_j}} = \int_Z \p{\int_{W_z} \chi v_j \frac{\partial f}{\partial x_j} \rho_z \mathrm{d}m_{W_z}} \mathrm{d}\nu(z),
\end{equation*}
where $\nu$ is a measure with finite mass on $Z$ and the $\rho_z$'s are uniformly $\mathcal{C}^{r-1}$ (and hence $\mathcal{C}^{1+q}$, see for instance \cite[Corollary 6.1.4]{ledrappierMetricEntropyDiffeomorphisms1985a}). Hence, we see that, for $z \in Z$, the function $\chi v_j \rho_z$ on $W_z$ is $\mathcal{C}^{1+q}$ and supported away from the boundary of $W_z$, so that we have by \cite[Lemma 3.2]{GLK}
\begin{equation*}
\va{\int_{W_z} \chi v_j \frac{\partial f}{\partial x_j} \rho_z \mathrm{d}m_{W_z}} \leq C \n{f}_{\gl},
\end{equation*}
for some constant $C > 0$ which is uniform in $z \in Z$ (because all data are). Integrating this estimate with respect to $\nu$, we find that
\begin{equation*}
\va{\mu\p{\chi v_j \frac{\partial f}{\partial x_j}}} \leq C \n{f}_{\gl}.
\end{equation*}
Summing over $j$, we find that $f  \mapsto \mu(\chi Vf)$ extends to a continuous linear from on $\mathcal{B}^{1,q}$, and then the lemma follows by a partition of unity argument.
\end{proof}

\begin{remark}
The two main elements in the proof of Lemma \ref{lemma:convergence} are Lemma \ref{lemma:linear_form} and the spectral decomposition \eqref{eq:spectral_decomposition} of the Koopman operator, valid on the space $\mathcal{B}^{1,q}$ (we recall here that we reversed time with respect to \cite{GLK} in order to work with the Koopman operator instead of the transfer operator). Such a decomposition for the Koopman operator has been established on many other spaces, see for instance \cite{Tsu0,Tsu,faureSemiclassicalApproachAnosov2008}. One can consequently wonder if Lemma \ref{lemma:linear_form} also holds on these spaces (it would give an alternative proof of Lemma \ref{lemma:convergence}). The answer seems to be positive for the spaces from \cite{Tsu0}, but the proof would be much more technically involved. The advantage of the space $\mathcal{B}^{1,q}$ from \cite{GLK} here is that it is easier to relate the linear form from Lemma \ref{lemma:linear_form} with the definition of these spaces.
\end{remark}

\section{Computing integrals with respect to $\mu$}\label{section:monte_carlo}

In order to achieve Step 4 of the S3 algorithm, we need to find a way to approximate the integral $\mu(g)$ of a Hölder observable $g$ against the SRB measure $\mu$ for $\phi$. We propose two simple methods to do so: a probabilistic one that relies on the physicality of $\mu$ (this is the method that is used in \cite{chandramoorthyEfficientComputationLinear2021}) and a deterministic one that uses a classical construction of the SRB measure. One could probably use more advanced techniques to compute integrals against $\mu$ (see for instance the method in \cite{crimminsFourierApproximationStatistical2020} to approximate $\mu$). 

In Remarks \ref{remark:precision_monte_carlo} and \ref{remark:precision_Riemann}, we discuss theoretical error estimates for the computation of the S3 algorithm. These estimates do not take into account, for instance, rounding errors due to finite precision arithmetic that would appear in an actual implementation of the algorithm. 

\subsection{Probabilistic method}\label{subsection:Montecarlo}
We describe here a ``Monte Carlo-like'' method to approximate integrals against $\mu$. This approach is quite standard (see for instance \cite{gottwaldSpuriousDetectionLinear2016} for another application of this method). In view of the property \eqref{eq:physical} of the measure $\mu$, it seems natural to approximate $\mu(g)$ by
\begin{equation}\label{eq:monte_carlo}
\mu(g) \simeq \frac{1}{N} \sum_{k = 0}^{N-1} g(\phi^k x),
\end{equation}
where $x$ is a point picked at random in $M$ and $N$ a large integer. This strategy works for Lebesgue almost every choice of $x$. The error in this approximation can be bounded almost surely using the following corollary of the almost sure invariance principle. We expect that non-asymptotic bound with high probability could be derived using a concentration inequality in the spirit of \cite{chazottesOptimalConcentrationInequalities2012}.

\begin{lemma}\label{lemma:monte_carlo}
Let $g$ be a Hölder-continuous function on $M$. Then, for Lebesgue almost every $x \in M$, there is a constant $C$ such that, for $N \geq 3$, we have
\begin{equation*}
\va{\frac{1}{N} \sum_{k = 0}^{N-1} g(\phi^k x) - \mu(g)} \leq C \sqrt{\frac{\log \log N}{N}}.
\end{equation*}
\end{lemma}

\begin{proof}
Replacing $g$ by $g - \mu(g)$, we can assume that $\mu(g) = 0$. If $M$ is endowed with any smooth measure, then the sequence of functions $(g, g \circ \phi, g \circ \phi^2,\dots)$ becomes a stochastic process and it follows from \cite[Theorem 2.1]{gouezelAlmostSureInvariance2010}, using the functional analytic tools from \cite{GLK,GouLiv2} or \cite{Tsu0,Tsu} as explained in \cite{gouezelAlmostSureInvariance2010}, that this stochastic process satisfies an almost sure principle. Consequently, there is a probability space $\Omega$, a process $(A_0,A_1,\dots)$ on $\Omega$ with the same law as $(g,g\circ \phi,\dots)$ and a Brownian motion $(B_t)_{t \geq 0}$ on $\Omega$ such that, almost surely,
\begin{equation*}
\va{\sum_{k = 0}^{N-1} A_k - B_{N}} \underset{N \to + \infty}{=} o(N^{\frac{1}{2}}).
\end{equation*}

From the law of the iterated logarithm, we know that almost surely there is a constant $C$ such that 
\begin{equation*}
\va{B_N} \leq C \sqrt{N \log \log N},
\end{equation*}
so that
\begin{equation*}
\va{ \frac{1}{N}\sum_{k = 0}^{N-1} A_k} \leq C' \sqrt{\frac{\log \log N}{N}},
\end{equation*}
for some constant $C' > C$ that depends on the alea. Since $(A_0,A_1,\dots)$ has the same law as $(g,g\circ \phi,\dots)$, the result follows.
\end{proof}

\begin{remark}\label{remark:precision_monte_carlo}
Let us go back to the S3 algorithm. If we use \eqref{eq:approximation_S3} as an approximation for $\Psi_{\phi}(X,f)$, the integrals being computed by the probabilistic method we just described, then we make three kinds of error:
\begin{itemize}
\item if we use orbits of length $n$ to compute approximate values of $\rho$ and $V$, then we make an error of size $\mathcal{O}(\eta^n)$, for some $\eta < 1$, when approximating them (see Remarks \ref{remark:solving_equations}, \ref{remark:step_2} and \ref{remark:compute_rho});
\item the truncation after the $L$th term in the series in the right hand side of \eqref{eq:approximation_S3} produces an error of size $\mathcal{O}(\eta^L)$;
\item the errors when approximating the integrals are almost surely of size $\mathcal{O}((\log \log N/N)^{1/2})$ by Lemma \ref{lemma:monte_carlo}.
\end{itemize}
Notice however that the constants in the $\mathcal{O}$ that appears in the third kind of error depend on the choice of $L$ (and on the alea) in an unexplicit way, so that this approach does not give an explicit bound on the error in terms of only the computation time.
\end{remark}

\subsection{Deterministic method}\label{subsection:deter}

The integral $\mu(g)$ may also be approximated by a deterministic method. Choose $W$ a $d_u$-dimensional submanifold of $M$ whose tangent space lieas in an unstable cone for $\phi$. For simplicity, we assume that $W$ belongs to $\Sigma$, the space of admissible unstable manifolds used to define the space $\mathcal{B}^{1,q}$ from \cite{GLK} in \S \ref{section:convergence}. Choose then a smooth function $\rho : W \to \R_+$ supported away from the boundary of $W$ and such that
\begin{equation*}
\int_W \rho \mathrm{d}m_W = 1.
\end{equation*}
Then, if $n$ is a very large integer, we can approximate $\mu(g)$ by
\begin{equation}\label{eq:approx_deter}
\mu(g) \simeq \int_W g \circ \phi^n \rho \mathrm{d}m_{W}.
\end{equation}
Indeed, the linear form
\begin{equation*}
u \mapsto \int_W u \rho \mathrm{d}m_W
\end{equation*}
is bounded on $\mathcal{B}^{1,q}$ (by the very definition of the space), and recalling the decomposition \eqref{eq:iterated_spectral_decomposition} for the Koopman operator we find that
\begin{equation*}
 \int_W g \circ \phi^n \rho \mathrm{d}m_{W} = \mu(g) +\int_W Q^n u \rho \mathrm{d}m_W.
\end{equation*}
Here, the last term decays exponentially fast with $n$ since the spectral radius of $Q$ is strictly less than $1$.

The integral $\int_W g \circ \phi^n \rho \mathrm{d}m_{W}$ can be approximated by Riemann sums. In the context of the S3 algorithm, the observable $g$ is only Hölder, so that it does not make sense to try to use a more sophisticated method to compute the integral. Moreover, when $n$ tends to $+ \infty$ the Hölder constant of $g \circ \phi^n$ grows exponentially fast, so that even the approximation by Riemann sum must be dealt with carefully. Let $RS_{N,n}(g)$ denote the approximation of $\int_W g \circ \phi^n \rho \mathrm{d}m_{W}$ by a Riemann sum with $N$ terms. If $g$ is $\alpha$-Hölder, then the restriction of $g \circ \phi^n$ to $W$ is $\alpha$-Hölder with a Hölder constant $\mathcal{O}(\kappa^{\alpha n})$, where $\kappa >1$ denotes the maximal expansion rate of $\phi$. Hence, we have
\begin{equation}\label{eq:quality_approximation}
\begin{split}
& \va{RS_{N,n}(g) - \mu(g)} \\  & \quad  \leq \va{RS_{N,n}(g) - \int_W g \circ \phi^n \rho \mathrm{d}m_{W}} + \va{\int_W g \circ \phi^n \rho \mathrm{d}m_{W} - \mu(g)} \\
    & \quad \lesssim N^{- \frac{\alpha}{d_u}} \kappa^{\alpha n} + \eta^n,
\end{split}
\end{equation}
with $\eta < 1$. We see that the best approximation here is obtained by taking $n \simeq \frac{\alpha \log N}{d_u \log(\kappa^\alpha/\eta)}$. With this value of $n$, the error when approximating $\mu(g)$ by $RS_{N,n}(g)$ is $\mathcal{O}(N^{- \beta})$ with $\beta = - \frac{\alpha \log \eta}{d_u \log(\kappa^\alpha/\eta)}$.

\begin{remark}\label{remark:precision_Riemann}
If $g$ is one of the observables for which we need to compute $\mu(g)$ in the fourth step of the S3 algorithm, then we assume that the value of $g$ at a point may be computed with an error $\mathcal{O}(\eta^n)$ in a time proportional to $n$ (one only needs to use an orbit of length proportional to $n$, see Remarks \ref{remark:solving_equations}, \ref{remark:step_2} and \ref{remark:compute_rho}). Hence, this error may be considered as part of the exponentially decaying terms in \eqref{eq:quality_approximation}. Computing $RS_{N,n}(g)$ requires the value of $g$ at $N$ points, hence the computation time $\tau_0$ of one such integral will be proportional to $n N = \mathcal{O}(N \log N)$. Consequently, up to machine precision, we compute one integral with an error $\mathcal{O}(\tau_0^{- \gamma})$, for some $\gamma > 0$, in a time $\tau_0$.

We also need to take into account the number of integrals that we want to compute: taking $L$ proportional to $\log \tau_0$ in \eqref{eq:unstable_contribution} gives an approximation of $\Psi_{\phi}(Y,f)$ with a precision $\mathcal{O}(\tau_0^{-\gamma})$. Here, the Hölder constants of the observables grow with $L$, but if the proportionality constant between $L$ and $\log \tau_0$ is small enough, this only makes the constant $\gamma$ smaller, by the analysis above. The total computation time $\tau$ is bounded by $\tau_0 \log \tau_0$. Up to taking $\gamma > 0$ smaller, the error in $\Psi_{\phi}(X,f)$ after time $\tau$ is $\mathcal{O}(\tau^{- \gamma})$, assuming that computations on orbits of length $n$ are made in a time proportional to $n$.
\end{remark}

\begin{remark}
Instead of using a piece of unstable manifold $W$, one could integrate directly against a smooth measure on the manifold $M$ (for instance if one is unable to find a piece of unstable manifold). The convergence also holds, but the quality of the approximation of an integral by Riemann sums deteriorates when the dimension gets bigger. Hence, we still get a polynomial precision, but \emph{a priori} with a worst exponent.
\end{remark}

\section{Extension to the case of hyperbolic attractors}\label{section:hyperbolic_attractors}

Let us explain now how our analysis of the S3 algorithm could be adapted to the case of \emph{hyperbolic attractors}, which happens to be the case considered in \cite{chandramoorthyEfficientComputationLinear2021}. Let $M$ be a Riemannian manifold, $U$ an open subset of $M$ and $\phi : U \to M$ a $\mathcal{C}^r$ embedding (we still assume $r > 2$). Let $ K \subseteq U$ be a compact, invariant, hyperbolic subset of $U$. We say that $K$ is a hyperbolic attractor for $\phi$ if, when the neighbourhood $U$ of $K$ is small enough, we have $K = \bigcap_{n \geq 0} \phi^n U$. In that case, we may also assume that $\phi(U) \subseteq U$, and we see then that for every $x \in U$ we have the distance between $\phi^n x$ and $K$ tends to $0$ when $n$ tends to $+ \infty$. Let us assume that $\phi_{|K}$ is transitive. Under these assumptions, $\phi$ admits a unique SRB measure $\mu$ supported in $K$ \cite{ruelle-attractor,Bow2} . The measure $\mu$ is physical, meaning that \eqref{eq:physical} holds for Lebesgue almost every $x \in U$ and every continuous function $f$ on $U$.

Moreover, Ruelle's formula \eqref{eq:Ruelle_formula} for the linear response is also satisfied in the context of hyperbolic attractors \cite{ruelleDifferentiationSRBStates1997,ruelleDifferentiationSRBStates2003,jiang_linear_response}. Consequently, the S3 algorithm may also be used in that case. Let us mention the few modifications that are needed to adapt the proof of Theorem \ref{theorem:main} in that case. 

In \S \ref{section:spectral_radii}, instead of considering sections of a vector bundle $E$ define on all $M$, we consider sections that are only defined on $K$. We can then define the space $\mathcal{C}^{\alpha,k}\p{K;E}$ as we did for $\mathcal{C}^{\alpha,k}\p{M;E}$ in the Anosov case, since for every $x \in K$ the unstable manifold of $x$ is contained in $K$. The proof of Proposition \ref{prop:bound_spectral_radius} adapts to this case without major changes.

Hence, if $X$ is a $\mathcal{C}^{1+}$ vector field defined on a neighbourhood of $K$, we may decompose $X$ as in \S \ref{section:Y_and_Z}
\begin{equation}\label{eq:decomposition_K}
X = Y + V - \phi_*V,
\end{equation}
with $Y$ and $V$ in $\mathcal{C}^{\alpha,1}\p{K;TM}$ and $Y$ parallel to the unstable direction. Notice however that in the case of a hyperbolic attractor the decomposition \eqref{eq:decomposition_K} only makes sense on the attractor $K$.

Since the SRB measure $\mu$ has absolutely continuous conditionals on the unstable manifold \cite[Theorem 1]{youngWhatAreSRB2002}, the analysis from \S \ref{section:unstable} still applies. In order to see that the study of the coboundary contribution from \S \ref{section:convergence} is still valid in the case of a hyperbolic attractor, one would just need to adapt the result from \cite{GLK} to the case of hyperbolic attractors. The technical tools from \cite{GouLiv2} can help (one could probably also use the spaces from \cite{Tsu0, Bal2}). Notice that in that case the Koopman operator \eqref{eq:koopman_operator} will be replaced by the operator
\begin{equation*}
\mathcal{K} : u \mapsto \chi u \circ \phi,
\end{equation*}
where $\chi$ is smooth, supported in $U$, and identically equal to $1$ on a neighbourhood of $K$. Then, the operator $P$ from \eqref{eq:spectral_decomposition} is not necessarily of rank $1$ anymore (unless $\phi$ is mixing), but one can ensure that the eigenvectors of $\mathcal{K}$ associated to eigenvalues of modulus $1$ are constant on a neighbourhood of $K$, so that $\mu(V P^n f) = 0$ for $f$ smooth and $n \in \mathbb{N}$, and the proof of Lemma \ref{lemma:convergence} still applies. To check that, use the spectral decomposition \cite[Theorem 18.3.1]{katokhasselblatt} to write $K = \bigcup_{j \in \mathbb{Z}/m\mathbb{Z}} K_j$ as a union of $m$ disjoint closed sets such that $\phi(K_j) = K_{j+1}$ and $\phi^m$ restricted to $K_j$ is mixing. Then, one can take $\chi$ of the form $\chi = \sum_{j = 1}^m \chi_j$ where $\chi_j$ takes value $1$ on a neighbourhood of $K_j$, the $\chi_j$ have disjoint supports, and the support of $\chi_{j+1} \circ \phi$ intersects the support of $\chi_\ell$ if and only if $\ell = j$. Then letting $\tilde{\chi}_j = \prod_{k = 0}^{+ \infty} \chi_{j + k} \circ \phi^{k}$, we can check that the eigenvectors of $\mathcal{K}$ associated to eigenvalues of modulus $1$ are spanned by the $\sum_{j = 1}^m e^{\frac{2 i j \ell \pi}{m}} \tilde{\chi}_j$ for $\ell \in \mathbb{Z} / m \mathbb{Z}$. This can be shown by reducing to the mixing case and using \cite[Theorem 7.5 and Lemma A.3]{Bal2}. Since $\tilde{\chi}_j$ is identically equal to $1$ near $K_j$ and to $0$ near $K_{\ell}$ for $\ell \neq j$, we have indeed that $\mu(V \tilde{\chi}_j) = 0$ for $V$ a vector field. The same kind of reasoning allows us to bypass the potential lack of exponential decay of correlations (in the absence of mixing) when evaluating the unstable contribution, so that the series \eqref{eq:unstable_contribution} is still exponentially converging. Notice that if $\phi$ is not mixing, and if one wants to use the deterministic method from \S \ref{subsection:deter} to compute integrals, then $\int_W g \circ \phi^n \rho\: \mathrm{d}m_{W}$ should be replaced by its Cesaro average (in $n$) in \eqref{eq:approx_deter}, which would give a much slower speed of convergence for the algorithm.

One of the main differences in the case of a hyperbolic attractor is when actually implementing: one may not be able to sample points in $K$ (which can be of zero Lebesgue measure). However, notice that in the actual implementation of the S3 algorithm, computations are made by using very long orbits for $\phi$. If $x$ is a point in $U$, then $\phi^n x$ converges to $K$ when $n$ tends to $+ \infty$ and the forward orbit of $x$ has the same asymptotic behavior as an orbit of a point in $K$. That is, there is $y \in K$ such that the distance between $\phi^n x$ and $\phi^n y$ tends to $0$ when $n$ tends to $+ \infty$ (see for instance \cite[Proposition 1.2 (c)]{ruelle-attractor}). Hence, for practical purposes, if one only uses very long orbits of points in $U$, it amounts to working with orbits of points in $K$ (since we are only working with Hölder functions, we can expect the error due to this approximation to be exponentially small in the length of the orbit we use, so that it should not harm the precision of the algorithm).

\bibliographystyle{alpha}
\bibliography{biblio.bib}

\newcommand{\etalchar}[1]{$^{#1}$}
\begin{thebibliography}{CFTW19}

\bibitem[AM07]{abramovMajda}
Rafail~V. Abramov and Andrew~J. Majda.
\newblock Blended response algorithms for linear fluctuation-dissipation for
  complex nonlinear dynamical systems.
\newblock {\em Nonlinearity}, 20(12):2793--2821, 2007.

\bibitem[Bal14]{baladiLinearResponseElse2014}
Viviane Baladi.
\newblock Linear response, or else.
\newblock In {\em Proceedings of the {{International Congress}} of
  {{Mathematicians}} ({{ICM}} 2014), {{Seoul}}, {{Korea}}, {{August}}
  13\textendash 21, 2014. {{Vol}}. {{III}}: {{Invited}} Lectures}, pages
  525--545. {KM Kyung Moon Sa}, {Seoul}, 2014.

\bibitem[Bal17]{quest}
Viviane Baladi.
\newblock The quest for the ultimate anisotropic {{Banach}} space.
\newblock {\em J. Stat. Phys.}, 166(3--4):525--557, 2017.

\bibitem[Bal18]{Bal2}
Viviane Baladi.
\newblock {\em Dynamical {{Zeta Functions}} and {{Dynamical Determinants}} for
  {{Hyperbolic Maps}}}, volume~68 of {\em Ergebnisse}.
\newblock {Springer}, 2018.

\bibitem[BGNN18]{bahsounRigorousComputationalApproach2017}
Wael Bahsoun, Stefano Galatolo, Isaia Nisoli, and Xiaolong Niu.
\newblock A rigorous computational approach to linear response.
\newblock {\em Nonlinearity}, 31(3):1073--1109, 2018.

\bibitem[BGW14]{blonigan2014least}
Patrick~J. Blonigan, Steven~A. Gomez, and Qiqi Wang.
\newblock Least squares shadowing for sensitivity analysis of turbulent fluid
  flows.
\newblock In {\em 52nd Aerospace Sciences Meeting}, page 1426, 2014.

\bibitem[BKL02]{BKL}
Michael Blank, Gerhard Keller, and Carlangelo Liverani.
\newblock Ruelle-{{Perron}}-{{Frobenius}} spectrum for {{Anosov}} maps.
\newblock {\em Nonlinearity}, 15(6):1905--1973, 2002.

\bibitem[BLL20]{bodai}
Tam{\'a}s B{\'o}dai, Valerio Lucarini, and Frank Lunkeit.
\newblock Can we use linear response theory to assess geoengineering
  strategies?
\newblock {\em Chaos: An Interdisciplinary Journal of Nonlinear Science},
  30(2):023124, 2020.

\bibitem[Bow70]{Bow2}
Rufus Bowen.
\newblock {\em Equilibrium {{States}} and the {{Ergodic Theory}} of {{Anosov
  Diffeomorphisms}}}.
\newblock Lecture {{Notes}} in {{Mathematics}}, 470. {Springer-Verlag,
  Berlin-Heidelberg-New York}, 1970.

\bibitem[BRS20]{bahsoun_ruziboev_saussol}
Wael Bahsoun, Marks Ruziboev, and Beno\^{\i}t Saussol.
\newblock Linear response for random dynamical systems.
\newblock {\em Adv. Math.}, 364:107011, 44, 2020.

\bibitem[BS12]{baladi-nonuniform}
Viviane Baladi and Daniel Smania.
\newblock Linear response for smooth deformations of generic nonuniformly
  hyperbolic unimodal maps.
\newblock {\em Annales scientifiques de l'ENS}, 45(6):861--926, 2012.

\bibitem[BS16]{bahsoun}
Wael Bahsoun and Benoît Saussol.
\newblock Linear response in the intermittent family: Differentiation in a
  weighted {$C^0$} -norm.
\newblock {\em Discrete \& Continuous Dynamical Systems}, 36(12):6657--6668,
  2016.

\bibitem[BT07]{Tsu0}
Viviane Baladi and Masato Tsujii.
\newblock Anisotropic {{H\"older}} and {{Sobolev}} spaces for hyperbolic
  diffeomorphisms.
\newblock {\em Ann. Inst. Fourier (Grenoble)}, 57(1):127--154, 2007.

\bibitem[BT08]{Tsu}
Viviane Baladi and Masato Tsujii.
\newblock Dynamical determinants and spectrum for hyperbolic diffemorphisms.
\newblock {\em Geometric and probabilistic structures in dynamics}, Amer. Math.
  Soc., Providence, RI(469):29--68, 2008.

\bibitem[BT16]{baladi-intermittent}
Viviane Baladi and Mike Todd.
\newblock Linear response for intermittent maps.
\newblock {\em Commun. Math. Phys.}, 347:857--–874, 2016.

\bibitem[Ces19]{cessac}
Bruno Cessac.
\newblock Linear response in neuronal networks: From neurons dynamics to
  collective response.
\newblock {\em Chaos}, 29:103105, 2019.

\bibitem[CF20]{crimminsFourierApproximationStatistical2020}
Harry Crimmins and Gary Froyland.
\newblock Fourier approximation of the statistical properties of {{Anosov}}
  maps on tori.
\newblock {\em Nonlinearity}, 33(11):6244--6296, 2020.

\bibitem[CFTW19]{ensembleSensitivity}
Nisha Chandramoorthy, Pablo Fernandez, Chaitanya Talnikar, and Qiqi Wang.
\newblock Feasibility analysis of ensemble sensitivity computation in turbulent
  flows.
\newblock {\em AIAA Journal}, 57(10):4514--4526, 2019.

\bibitem[CG12]{chazottesOptimalConcentrationInequalities2012}
Jean-Ren{\'e} Chazottes and S{\'e}bastien Gou{\"e}zel.
\newblock Optimal {{Concentration Inequalities}} for {{Dynamical Systems}}.
\newblock {\em Communications in Mathematical Physics}, 316(3):843--889, 2012.

\bibitem[CL99]{chiconeEvolutionSemigroupsDynamical1999}
Carmen Chicone and Yuri Latushkin.
\newblock {\em Evolution {{Semigroups}} in {{Dynamical Systems}} and
  {{Differential Equations}}}, volume~70 of {\em Mathematical {{Surveys}} and
  {{Monographs}}}.
\newblock {American Mathematical Society}, {Providence, Rhode Island}, 1999.

\bibitem[CW20]{chandramoorthyComputableRealizationRuelle2020}
Nisha Chandramoorthy and Qiqi Wang.
\newblock A computable realization of {{Ruelle}}'s formula for linear response
  of statistics in chaotic systems.
\newblock {\em arXiv:2002.04117}, 2020.

\bibitem[CW21]{chandramoorthyEfficientComputationLinear2021}
Nisha Chandramoorthy and Qiqi Wang.
\newblock Efficient computation of linear response of chaotic attractors with
  one-dimensional unstable manifolds.
\newblock {\em (to appear in SIAM Journal on Applied Dynamical Systems)},
  arXiv:2103.08816, 2021.

\bibitem[Den89]{denkerCentralLimitTheorem1989}
Manfred Denker.
\newblock {\em The Central Limit Theorem for Dynamical Systems}.
\newblock 1989.
\newblock Dynamical systems and ergodic theory, 28th Sem. St. Banach Int. Math.
  Cent., Warsaw/Pol. 1986, Banach Cent. Publ. 23, 33-62 (1989).

\bibitem[DLL01]{delallaveRemarksSobolevRegularity2001}
Rafael De~La~Llave.
\newblock Remarks on {{Sobolev}} regularity in {{Anosov}} systems.
\newblock {\em Ergodic Theory and Dynamical Systems}, 21(04), 2001.

\bibitem[Dol04]{dolgopyat}
Dmitry Dolgopyat.
\newblock On differentiability of {SRB} states for partially hyperbolic
  systems.
\newblock {\em Inventiones Mathematicae}, 155(2):389--449, 2004.

\bibitem[DW11]{dow2011uncertainty}
Eric Dow and Qiqi Wang.
\newblock Uncertainty quantification of structural uncertainties in rans
  simulations of complex flows.
\newblock In {\em 20th AIAA Computational Fluid Dynamics Conference}, page
  3865, 2011.

\bibitem[EHL04]{eyink}
Gregory Eyink, Thomas Haine, and Daniel~J. Lea.
\newblock Ruelle{\textquotesingle}s linear response formula, ensemble adjoint
  schemes and {L}{\'{e}}vy flights.
\newblock {\em Nonlinearity}, 17(5):1867--1889, 2004.

\bibitem[FD11]{fidkowski2011review}
Krzysztof~J. Fidkowski and David~L. Darmofal.
\newblock Review of output-based error estimation and mesh adaptation in
  computational fluid dynamics.
\newblock {\em AIAA journal}, 49(4):673--694, 2011.

\bibitem[FRS08]{faureSemiclassicalApproachAnosov2008}
Fr{\'e}d{\'e}ric Faure, Nicolas Roy, and Johannes Sj{\"o}strand.
\newblock Semi-classical approach for {{Anosov}} diffeomorphisms and {{Ruelle}}
  resonances.
\newblock {\em Open Mathematics Journal}, 1:35--81, 2008.

\bibitem[GG19]{galatoloLinearResponseDynamical2017}
Stefano Galatolo and Paolo Giulietti.
\newblock A linear response for dynamical systems with additive noise.
\newblock {\em Nonlinearity}, 32:2269--2301, 2019.

\bibitem[GHL04]{ghl_geometry}
Sylvestre Gallot, Dominique Hulin, and Jacques Lafontaine.
\newblock {\em Riemannian geometry}.
\newblock Universitext. Springer-Verlag, Berlin, third edition, 2004.

\bibitem[GL06]{GLK}
Sebastien Gou{\"e}zel and Carlangelo Liverani.
\newblock Banach spaces adapted to {{Anosov}} systems.
\newblock {\em Ergod. Th. and Dyn. Sys.}, 26:189--217, 2006.

\bibitem[GL08]{GouLiv2}
S{\'e}bastien Gou{\"e}zel and Carlangelo Liverani.
\newblock Compact locally maximal hyperbolic sets for smooth maps: Fine
  statistical properties.
\newblock {\em J. Differential Geom.}, 79(3):433--477, 2008.

\bibitem[GL20]{ghilLucarini}
Michael Ghil and Valerio Lucarini.
\newblock The physics of climate variability and climate change.
\newblock {\em Rev. Mod. Phys.}, 92:035002, Jul 2020.

\bibitem[Gou10]{gouezelAlmostSureInvariance2010}
S{\'e}bastien Gou{\"e}zel.
\newblock Almost sure invariance principle for dynamical systems by spectral
  methods.
\newblock {\em The Annals of Probability}, 38(4):1639--1671, 2010.

\bibitem[GWW16]{gottwaldSpuriousDetectionLinear2016}
Georg~A. Gottwald, Caroline~L. Wormell, and Jeroen Wouters.
\newblock On spurious detection of linear response and misuse of the
  fluctuation\textendash dissipation theorem in finite time series.
\newblock {\em Physica D: Nonlinear Phenomena}, 331:89--101, 2016.

\bibitem[HM10]{hairer}
Martin Hairer and Andrew~J. Majda.
\newblock A simple framework to justify linear response theory.
\newblock {\em Nonlinearity}, 23(4):909--922, 2010.

\bibitem[HM20]{huhn}
Francisco Huhn and Luca Magri.
\newblock Stability, sensitivity and optimisation of chaotic acoustic
  oscillations.
\newblock {\em Journal of Fluid Mechanics}, 882:A24, 2020.

\bibitem[Hus94]{husemoller}
Dale Husemoller.
\newblock {\em Fibre bundles}, volume~20 of {\em Graduate Texts in
  Mathematics}.
\newblock Springer-Verlag, New York, third edition, 1994.

\bibitem[Jia12]{jiang_linear_response}
Miaohua Jiang.
\newblock Differentiating potential functions of {SRB} measures on hyperbolic
  attractors.
\newblock {\em Ergodic Theory Dynam. Systems}, 32(4):1350--1369, 2012.

\bibitem[KH95]{katokhasselblatt}
Anatole Katok and Boris Hasselblatt.
\newblock {\em Introduction to the Modern Theory of Dynamical Systems},
  volume~54 of {\em Encyclopedia of {{Mathematics}} and Its {{Applications}}}.
\newblock {Cambridge University Press, Cambridge}, 1995.

\bibitem[KL99]{Ke-Li}
Gerhard Keller and Carlangelo Liverani.
\newblock Stability of the spectrum for transfer operators.
\newblock {\em Ann. Scuola Norm. Sup. Pisa Cl. Sci.}, 28(1):141--152, 1999.

\bibitem[LBH{\etalchar{+}}14]{lucarini1}
Valerio Lucarini, Richard Blender, Corentin Herbert, Salvatore Pascale,
  Francesco Ragone, and Jeroen Wouters.
\newblock Mathematical and physical ideas for climate science.
\newblock {\em Rev. Geophys.}, 52:809--859, 2014.

\bibitem[LY85]{ledrappierMetricEntropyDiffeomorphisms1985a}
Francois Ledrappier and Lai-Sang Young.
\newblock The {{Metric Entropy}} of {{Diffeomorphisms}}: {{Part I}}:
  {{Characterization}} of {{Measures Satisfying Pesin}}'s {{Entropy Formula}}.
\newblock {\em The Annals of Mathematics}, 122(3):509, 1985.

\bibitem[Ni19]{angxiuFluid}
Angxiu Ni.
\newblock Hyperbolicity, shadowing directions and sensitivity analysis of a
  turbulent three-dimensional flow.
\newblock {\em Journal of Fluid Mechanics}, 863:644–669, 2019.

\bibitem[Ni21]{ni2021fast}
Angxiu Ni.
\newblock Fast linear response algorithm for differentiating chaos.
\newblock {\em arXiv:2009.00595}, 2021.

\bibitem[NW17]{angxiuNILSS}
Angxiu Ni and Qiqi Wang.
\newblock Sensitivity analysis on chaotic dynamical systems by non-intrusive
  least squares shadowing (nilss).
\newblock {\em Journal of Computational Physics}, 347:56--77, 2017.

\bibitem[OP89]{oreyDeviationsTrajectoryAverages2021}
Steven {Orey} and Stephan {Pelikan}.
\newblock {Deviations of trajectory averages and the defect in Pesin's formula
  for Anosov diffeomorphisms}.
\newblock {\em {Trans. Am. Math. Soc.}}, 315(2):741--753, 1989.

\bibitem[RLL16]{lucarini2}
Francesco Ragone, Valerio Lucarini, and Frank Lunkeit.
\newblock A new framework for climate sensitivity and prediction: a modelling
  perspective.
\newblock {\em Climate Dynamics}, 46:1459--1471, 2016.

\bibitem[Rue76]{ruelle-attractor}
David Ruelle.
\newblock A measure associated with axiom-{A} attractors.
\newblock {\em Amer. J. Math.}, 98(3):619--654, 1976.

\bibitem[Rue97]{ruelleDifferentiationSRBStates1997}
David Ruelle.
\newblock Differentiation of {{SRB States}}.
\newblock {\em Communications in Mathematical Physics}, 187(1):227--241, 1997.

\bibitem[Rue03]{ruelleDifferentiationSRBStates2003}
David Ruelle.
\newblock Differentiation of {{SRB}} states: {{Correction}} and complements.
\newblock {\em Communications in Mathematical Physics}, 234(1):185--190, 2003.

\bibitem[SF18]{shimizu2018output}
Yukiko~S. Shimizu and Krzysztof Fidkowski.
\newblock Output-based error estimation for chaotic flows using reduced-order
  modeling.
\newblock In {\em 2018 AIAA Aerospace Sciences Meeting}, page 0826, 2018.

\bibitem[Sin68]{sinaiMarkovPartitionsDiffeomorphisms1968}
Yakov~G. Sinai.
\newblock Markov partitions and {C} -diffeomorphisms.
\newblock {\em Functional Analysis and its Applications}, 2:61--82, 1968.

\bibitem[{Sin}73]{sinaiGIBBSMEASURESERGODIC}
Yakov~G. {Sinai}.
\newblock {Gibbs measures in ergodic theory}.
\newblock {\em {Russ. Math. Surv.}}, 27(4):21--69, 1973.

\bibitem[SW21]{sliwiak2021spacesplit}
Adam~A. Sliwiak and Qiqi Wang.
\newblock Space-split algorithm for sensitivity analysis of discrete chaotic
  systems with unstable manifolds of arbitrary dimension.
\newblock {\em arXiv:2109.13313}, 2021.

\bibitem[WG18]{wormell2}
Caroline~L. Wormell and Georg~A. Gottwald.
\newblock On the validity of linear response theory in high-dimensional
  deterministic dynamical systems.
\newblock {\em J Stat Phys}, 172:1479–--1498, 2018.

\bibitem[WG19]{wormell1}
Caroline~L. Wormell and Georg~A. Gottwald.
\newblock Linear response for macroscopic observables in high-dimensional
  systems.
\newblock {\em Chaos: An Interdisciplinary Journal of Nonlinear Science},
  29(11):113127, 2019.

\bibitem[Yoc95]{yoccozIntroductionHyperbolicDynamics1995}
Jean-Christophe Yoccoz.
\newblock Introduction to hyperbolic dynamics.
\newblock In {\em Real and Complex Dynamical Systems. {{Proceedings}} of the
  {{NATO Advanced Study Institute}} Held in {{Hiller\o d}}, {{Denmark}},
  {{June}} 20-{{July}} 2, 1993}, pages 265--291. {Kluwer Academic Publishers},
  {Dordrecht}, 1995.

\bibitem[You02]{youngWhatAreSRB2002}
Lai-Sang Young.
\newblock What are {{SRB}} measures, and which dynamical systems have them?
\newblock In {\em Journal of {{Statistical Physics}}}, volume 108, pages
  733--754. 2002.

\end{thebibliography}

\end{document}